\documentclass{amsart}
\usepackage{amssymb}
\usepackage{graphicx}
\usepackage{latexsym,amscd}
\usepackage{amsmath,paralist}%,diagrams}
\usepackage{t1enc}
\title[Narrow  2-Class Field Towers]
          { On the Narrow $2$-Class Field Tower of  \\ Some Real Quadratic Number Fields \\ Part I: Ranks}
\author{E.\,Benjamin,  C.\,Snyder}

\keywords{real quadratic field, Hilbert class field, 2-class group, narrow class group, narrow 2-class field tower}
\subjclass[2010]{11R29, 20D15}

\begin{document}
\newcommand{\rsp}{\raisebox{0em}[2.7ex][1.3ex]{\rule{0em}{2ex} }}
\newcommand{\ul}{\underline}
\newcommand{\1}{\boldsymbol{1}}
\newcommand{\I}{\boldsymbol{i}}
\newcommand{\J}{\boldsymbol{j}}
\newcommand{\K}{\boldsymbol{k}}
\newcommand{\C}{{\mathbb C}}
\newcommand{\h}{\mathbb{H}}
\newcommand{\Q}{{\mathbb Q}}
\newcommand{\R}{{\mathbb R}}
\newcommand{\B}{{\mathbb B}}
\newcommand{\F}{{\mathbb F}}
\newcommand{\N}{{\mathbb N}}
\newcommand{\Z}{{\mathbb Z}}
\newcommand{\aQ}{{\overline{\Q}}}
\newcommand{\cA}{\mathcal A}
\newcommand{\OO}{\mathcal O}
\newcommand{\Rl}{{\rm Re}}
\newcommand{\cM}{\mathcal M}
\newcommand{\Am}{\operatorname{Am}}
\newcommand{\eq}{\stackrel{2}{=}}
\newcommand{\fZ}{\mathfrak 2}
\newcommand{\fQ}{\mathfrak Q}
\newcommand{\fa}{\mathfrak a}
\newcommand{\fb}{\mathfrak b}
\newcommand{\fc}{\mathfrak c}
\newcommand{\fC}{\mathfrak C}
\newcommand{\fE}{\mathfrak E}
\newcommand{\fp}{\mathfrak p}
\newcommand{\fP}{\mathfrak P}
\newcommand{\fA}{\mathfrak A}
\newcommand{\fB}{\mathfrak B}
\newcommand{\fG}{\mathfrak G}
\newcommand{\fr}{\mathfrak r}
\newcommand{\frb}{\widetilde{\fr}}
\newcommand{\fR}{\mathfrak R}
\newcommand{\fq}{\mathfrak q}
\newcommand{\gac}{\widehat}
\newcommand{\Ft}{\widetilde{F}}
\newcommand{\disc}{\operatorname{disc}}
\newcommand{\ch}{\operatorname{char}}
\newcommand{\id}{\operatorname{id}}
\newcommand{\rank}{\operatorname{rank}}
\newcommand{\lcm}{{\operatorname{lcm}}}
\newcommand{\gen}{{\operatorname{gen}}}
\newcommand{\cen}{{\operatorname{cen}}}
\newcommand{\sfk}{{\operatorname{sfk}}}
\newcommand{\Tr}{\operatorname{Tr}}
\newcommand{\Cl}{\operatorname{Cl}}
\newcommand{\Ram}{\operatorname{Ram}}
\newcommand{\Gal}{\operatorname{Gal}}
\newcommand{\Aut}{\operatorname{Aut}}
\newcommand{\AGL}{\operatorname{AGL}}
\newcommand{\Ind}{\operatorname{Ind}}
\newcommand{\GL}{\operatorname{GL}}
\newcommand{\im}{\operatorname{im}\,}
\newcommand{\eps}{\varepsilon}
\newcommand{\la}{\langle}
\newcommand{\ra}{\rangle}
\newcommand{\lra}{\rightarrow}
\newcommand{\too}{\mapsto}
\newcommand{\ts}{\textstyle}
\newcommand{\ab}{\operatorname{ab}}
\newcommand{\ff}{\mathfrak f}
\newcommand{\fm}{\mathfrak m}
\newcommand{\fn}{\mathfrak n}
\def\la{\langle}
\def\ra{\rangle}
\def\lra{\rightarrow}
\def\lms{\mapsto}
\def\vp{\varphi}
\def\vt{\vartheta}
\def\tr{\mbox{{\rm Tr}}\,}
\def\hom{\mbox{{\rm Hom}}\,}
\def\ind{\mbox{{\rm Ind}}}
\newtheorem{df}{Definition}
\newtheorem{thm}{Theorem}
\newtheorem{lem}{Lemma}
\newtheorem{prop}{Proposition}
\newtheorem{cor}{Corollary}

\begin{abstract} We determine some properties of the narrow 2-class field tower of those real quadratic number fields  whose  discriminants are not a sum of two squares and for which their 2-class groups are elementary of order $4$. Here in Part I, we determine precisely  which of these  fields have the  2-class groups of their narrow 2-Hilbert class fields of rank 2.
\end{abstract}

\maketitle

\section{{\bf Introduction}}

Let $k$ be an algebraic  number field and consider the narrow 2-class field tower over $k$. This is the sequence of fields $\{k^n_+\},$ where $k^0_+=k$ and $k^{n+1}_+$ is the narrow 2-class field extension of $k_+^n$, i.e., the maximal abelian 2-extension of $k^n_+$ unramified outside $\infty$  (hence unramified at all the finite primes of $k^n_+$).
When $k$ is, for instance, an  imaginary quadratic field, the narrow 2-class field tower coincides with the ordinary  2-class field tower $\{k^n\}$.
However,  when $k$ is real,  then of course the  two towers can differ quite radically. We are now interested in continuing our comparison of  these two towers for the family of real quadratic fields $k$  for which the ordinary  $2$-class group $\Cl_2(k)\simeq (2,2)=V_4$, the Klein $4$-group. The reason for considering this family is that the ordinary $2$-class field tower for any quadratic field $k$ with $\Cl_2(k)\simeq (2,2)$ is  completely understood; see \cite{Kis}, \cite{CD}, and \cite{BS}. Now when the discriminant $d_k$ of $k$ is a sum of two squares (equivalently, $d_k$ is the product of  positive prime discriminants), then the ordinary, and in some cases the narrow towers, were studied in  \cite{CD} and then extended to all the narrow towers in \cite{BS1}.

This brings us to the purpose of the present paper. Ultimately, we are interested in studying the length  of  the narrow $2$-class field tower of those real quadratic number fields $k$  for which $\Cl_2(k)\simeq (2,2)$ and for which the discriminant $d_k$ is not the sum of two squares. More precisely, we want to determine  which $k$ have narrow $2$-class field towers of length $2$, in which case the rest have tower length $\geq 3$. As we will see,  $\rank\Cl_2(k_+^1)=2$ implies that the length of the tower is $2$, whereas if the rank is $\geq 3$, the tower length can be $2$ or $\geq 3$. The main object of Part I here is to determine precisely which of these fields $k$ have $\rank\Cl_2(k_+^1)=2$ and hence determine those fields  which  ``trivially'' have finite $2$-class field towers.

Namely, we prove the following main theorem. Let $G^+=\Gal(k_+^2/k)$ and $G^+_3$ the third term in the lower central series of $G^+$ (see Section~\ref{S6} below for details).

\bigskip
\noindent{\bf Main Theorem.} {\em Let $k$ be a real quadratic number field for which $\Cl_2(k)\simeq (2,2)$ and for which the discriminant $d_k$ is not the sum of two squares. Then
$$\rank\Cl_2(k_+^1)\left\{
                     \begin{array}{ll}
                       \geq 3, &  \hbox{if} \;\;(G^+:G_3^+)=64,\\
                       =3, & \hbox{if} \;\;G^+/G_3^+\simeq 32.033,\\
                       =2, & \hbox{otherwise},
                     \end{array}
                   \right.$$
where $32.033$ denotes the group as coded in \cite{HS}.}

\bigskip
In Part II we will then determine (except for possibly one family of fields) exactly which fields have narrow $2$-tower length $\geq 3$. Hence  in these cases we have the possibility of infinite narrow $2$-class field tower.

Now back to Part I. We start by exhibiting all those real quadratic number fields, $k$, for which $d_k$ is not the sum of two squares and such that $\Cl_2(k)$ has rank $2$. These fields are then partitioned into their $4$-ranks: $2,1,0$. We then concentrate on the elementary case, i.e., when $4$-$\rank\Cl_2(k)=0$.

\section{{\bf Real Quadratic Number Fields with ${2}$-Class Groups of Rank ${2}$ and whose  Discriminants are Not a Sum of Two Squares  }}\label{S2}

We now consider real quadratic number fields $k$ with $\rank\Cl_2(k)=2$ and the discriminant $d_k$ not a sum of two squares.  Then by genus theory $d_k$ is a product of four prime discriminants $d_j$: $d_k=d_1d_2d_3d_4$, for which either precisely two or all four are negative.  We partition the fields into four types according  to the number of negative prime discriminants dividing  $d_k$ and whether or not $-4$ is one of them:

\bigskip
\noindent Type I: $d_1,d_2>0$; $d_3,d_4<0$; $d_j\not=-4$, ($j=3,4$);

\noindent Type II: $d_1,d_2>0$; $d_3,d_4<0$; $d_4=-4$;

\noindent Type III: $d_j<0$, $d_j\not=-4$, ($j=1,2,3,4$);

\noindent Type IV: $d_j<0$, $d_4=-4$, ($j=1,2,3,4$).

\bigskip

We are interested in the cases where $\Cl_2(k)$ is elementary (hence $\Cl_2(k)\simeq (2,2)$ and so the narrow $2$-class group $\Cl_2^+(k)\simeq (2,2,2)$). We start more generally by determining the $4$-rank of  $\Cl_2(k)$. As in some of our previous work, we use diagrams on the prime discriminants dividing $d_k$ to give  information that we wish to display, cf., e.g., \cite{BLS2}.  Let $p_i$ denote the prime dividing $d_i$, for $i=1,2,3,4$. Then we exhibit an arrow from $d_i$ to $d_j$ iff $(d_i/p_j)=-1$, where $(./.)$ denotes the Kronecker symbol. If there is an arrow in both directions, we simply exhibit a line segment between the two prime discriminants.

Appendix I shows a complete list of the discriminants of all real  quadratic number fields $k$ with $\rank\Cl_2(k)=2$ and where $d_k$ is not a sum of two squares.  The discriminants in this list are partitioned first by type, as given above, and then within each type, by the $4$-rank of $\Cl_2(k)$.   The four points in each diagram represent the prime discriminants $d_1,d_2,d_3,d_4$, ordered as in the list of types above, such that $d_1$ is the upper left point, $d_2$ upper right, $d_3$ lower left,  and $d_4$ lower right. (For those readers who have no interest in the diagrams, we list in Appendix I$'$ the values of the relevant Kronecker symbols.)

To compute the $4$-rank $r$ of $\Cl_2(k)$, (which in our case is equal to the $4$-rank of $\Cl_2^+(k)$, since $\Cl_2^+(k)\simeq C_2\times \Cl_2(k)$)  recall that $2^r$ is the number of $C_4$-splittings of $d_k$; see \cite{RR} or \S 2 of \cite{Lem}.  A factorization $d_k=D_1D_2$ into two relatively prime discriminants (including $d_k=1\cdot d_k$) is a $C_4$-splitting, if $(D_1/p_2)=(D_2/p_1)=1$ for all $p_i\mid D_i$.  As an example of this, in Appendix I consider Type I($\alpha$).
It is convenient using the diagram to see that there are four $C_4$-splittings of $d_k$: namely,
$$\{D_1,D_2\}\in\big\{\{1,d_k\},\{d_1,d_2d_3d_4\},\{d_2,d_1d_3d_4\},\{d_1d_2,d_3d_4\}\big\},$$
which implies that  $4$-$\rank\Cl_2(k)=2$ in this case.

We now show that the lists in Appendix I  are complete. (This will show  another advantage of using the diagrams.)

We first consider  Type I. There are $2^6=64$ possible diagrams depending on the number of line segments occurring between the six pairs of points, but we may without loss of generality assume that $(d_4/p_3)=-1$, in which case we are reduced to $32$ diagrams.  We only exhibit inequivalent graphs, where one diagram is equivalent to another if there is a bijection of the set of the points  $\{d_1,d_2\}$ to itself   taking one diagram to the other. Here is a list of the cases (as labeled in the list in Appendix I) and number of equivalent diagrams for each case, the sum of which is $32$:

\bigskip

\bigskip
\begin{tabular}{l|c|cccccc|ccc}
Type I & $(\alpha)$ & ($A_1$) & ($A_2$) & ($A_3$) & ($A_4$) & ($A_5$) & ($A_6$) & ($a_1$) & ($a_2$) & ($a_3$)  \\ \hline
Number &    1  & 1 & 2 & 2 & 2 & 1 & 1 & 1 & 1 & 2 \\
\end{tabular}

\bigskip
\begin{tabular}{l|cccccccccc}
Type I & ($a_4$) & ($a_5$) & ($a_6$) & ($a_7$) & ($a_8$) & ($a_9$) & ($a_{10}$) & ($a_{11}$) & ($a_{12}$) & ($a_{13}$) \\ \hline
Number &  2 &   2 &  2 & 2 & 2 & 1 & 1 &  2 & 2 & 2  \\
\end{tabular}

\bigskip
\bigskip

Of the $64$  Type II  diagrams, we combine each  pair for which $(d_3/2) (=(d_3/p_4))= \pm 1$. Moreover, since $(-4/p_3)=-1$, we depict the bottom part of each pair of graphs as:
\begin{picture}(10,10)(-5,-3)
    \thicklines
    %\put(20,20){\circle*{3}}
    %\put(0,20){\circle*{3}}
    \put(0,0){\circle*{3}}
    \put(20,0){\circle*{3}}
    %\put(0,20){\line(1,0){20}}
    \put(0,0){\line(1,0){20}}
    \put(0,-3){\line(1,0){20}}
    %\put(0,0){\line(0,1){20}}
    %\put(20,0){\line(-1,1){20}}
     %\put(0,0){\line(1,1){20}}
     %\put(20,0){\line(0,1){20}}
     \put(5,-3){\vector(-1,0){0}}
  \end{picture}
\hspace{2em} We will see below that the data we seek is (essentially) identical for  the graphs in each of these pairs.

\medskip
(Also notice that $(d_j/2)=\pm 1$, whereas $(-4/p_j)=+1$, for $j=1,2$.) Of the $32$ graphs, we have given a list of those which are inequivalent (equivalence as in the Type I above) in  Appendix I. Once again  we give a list of the cases  and number of equivalent diagrams for each case, the sum of which is $32$:

\bigskip

\bigskip
\begin{tabular}{l|c|cccccc|ccc}
Type II & $(\beta)$ & ($B_1$) & ($B_2$) & ($B_3$) & ($B_4$) & ($B_5$) & ($B_6$) & ($b_1$) & ($b_2$) & ($b_3$)  \\ \hline
Number &    1  & 1 & 2 & 2 & 2 & 1 & 1 & 1 & 1 & 2 \\
\end{tabular}

\bigskip
\begin{tabular}{l|cccccccccc}
Type II & ($b_4$) & ($b_5$) & ($b_6$) & ($b_7$) & ($b_8$) & ($b_9$) & ($b_{10}$) & ($b_{11}$) & ($b_{12}$) & ($b_{13}$) \\ \hline
Number &  2 &   2 &  2 & 2 & 2 & 1 & 1 &  2 & 2 & 2  \\
\end{tabular}

\bigskip

\bigskip
\bigskip

Now consider Type III. We may assume without loss of generality that all the graphs contain the subgraph

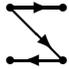
\begin{picture}(40,40)(-5,-8)
    \thicklines
    \put(20,20){\circle*{3}}
    \put(0,20){\circle*{3}}
    \put(0,0){\circle*{3}}
    \put(20,0){\circle*{3}}
    \put(0,20){\line(1,0){20}}
    %\put(0,0){\line(0,1){20}}
    \put(20,0){\line(-1,1){20}}
    %\put(0,0){\line(1,1){20}}
    %\put(20,0){\line(0,1){20}}
    \put(0,0){\line(1,0){20}}
    %\put(12,8){\vector(-1,1){0}}
   % \put(20,16){\vector(0,1){0}}
    %\put(4,4){\vector(-1,-1){0}}
    \put(15,20){\vector(1,0){0}}
    %\put(0,8){\vector(0,-1){0}}
    \put(16,4){\vector(1,-1){0}}
    \put(5,0){\vector(-1,0){0}}
  \end{picture}
For first we may, if necessary, switch the roles of $d_1$ with $d_2$ and $d_3$ with $d_4$ to insure that the horizontal arrows on top and bottom are oriented in the desired direction. If the diagonal arrow is headed upward, then rotate the diagram $180$ degrees to obtain the desired orientation of the subgraph.
Hence there are eight possible graphs depending on the direction of the arrows on the three remaining pairs of points, but not all are inequivalent. This time equivalence is given as above with respect to bijections on the set of all four prime discriminants. Here is a list of the cases and number of equivalent diagrams for each case, the sum of which is $8$:

\bigskip

\bigskip
\begin{tabular}{l|c|ccc}
Type III & ($C$) & ($c_1$) & ($c_2$) & ($c_3$)   \\ \hline
Number &    3  & 3 & 1 & 1  \\
\end{tabular}
\bigskip

Since for ($C$) and ($c_1$) it may not be clear as to which of the eight graphs are equivalent and why, we give them explicitly here:  $(C)\sim (C')\sim (C'')$  (see the diagrams below) with $\sigma'=(2\,4)$ and $\sigma''=(2\,3)$ where $d_j\mapsto d_{\sigma'(j)}$ gives the bijection of the prime discriminants, which transforms $(C)$ to $(C')$ and similarly for $\sigma''$ transforming $(C)$ to $(C'')$. Moreover $(c_1)\sim (c_1')\sim (c_1'')$  (again see the diagrams below) with $\tau'=(2\,4)$ and $\tau''=(1\,3\,4)$ (with the  $\tau$'s acting on the indices of the $d_j$).

\bigskip

Here are the diagrams:

\bigskip

\noindent \begin{tabular}{lll}\vspace{-.1in}
($C$)  &   ($C'$)  &   ($C''$)\\
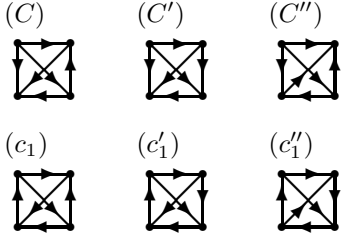
\begin{picture}(40,40)(-5,-8)  %graph for III(C)
    \thicklines
    \put(20,20){\circle*{3}}
    \put(0,20){\circle*{3}}
    \put(0,0){\circle*{3}}
    \put(20,0){\circle*{3}}
    \put(0,20){\line(1,0){20}}
    \put(0,0){\line(0,1){20}}
    \put(20,0){\line(-1,1){20}}
    \put(0,0){\line(1,1){20}}
    \put(20,0){\line(0,1){20}}
    \put(0,0){\line(1,0){20}}
    %\put(12,8){\vector(-1,1){0}}
    \put(20,16){\vector(0,1){0}}
    \put(4,4){\vector(-1,-1){0}}
    \put(15,20){\vector(1,0){0}}
    \put(0,8){\vector(0,-1){0}}
    \put(16,4){\vector(1,-1){0}}
    \put(5,0){\vector(-1,0){0}}
  \end{picture} &
  \begin{picture}(40,40)(-5,-8)  %graph for III(C')
    \thicklines
    \put(20,20){\circle*{3}}
    \put(0,20){\circle*{3}}
    \put(0,0){\circle*{3}}
    \put(20,0){\circle*{3}}
    \put(0,20){\line(1,0){20}}
    \put(0,0){\line(0,1){20}}
    \put(20,0){\line(-1,1){20}}
    \put(0,0){\line(1,1){20}}
    \put(20,0){\line(0,1){20}}
    \put(0,0){\line(1,0){20}}
    %\put(12,8){\vector(-1,1){0}}
    \put(20,8){\vector(0,-1){0}}
    \put(4,4){\vector(-1,-1){0}}
    \put(15,20){\vector(1,0){0}}
    \put(0,8){\vector(0,-1){0}}
    \put(16,4){\vector(1,-1){0}}
    \put(5,0){\vector(-1,0){0}}
  \end{picture} &
  \begin{picture}(40,40)(-5,-8)  %graph for III(C'')
    \thicklines
    \put(20,20){\circle*{3}}
    \put(0,20){\circle*{3}}
    \put(0,0){\circle*{3}}
    \put(20,0){\circle*{3}}
    \put(0,20){\line(1,0){20}}
    \put(0,0){\line(0,1){20}}
    \put(20,0){\line(-1,1){20}}
    \put(0,0){\line(1,1){20}}
    \put(20,0){\line(0,1){20}}
    \put(0,0){\line(1,0){20}}
    %\put(12,8){\vector(-1,1){0}}
    \put(20,16){\vector(0,1){0}}
    \put(9,9){\vector(1,1){0}}
    \put(15,20){\vector(1,0){0}}
    \put(0,8){\vector(0,-1){0}}
    \put(16,4){\vector(1,-1){0}}
    \put(5,0){\vector(-1,0){0}}
  \end{picture}\\
 \end{tabular}

 \noindent \begin{tabular}{lll}\vspace{-.1in}
($c_1$)  &   ($c_1'$)  &   ($c_1''$)\\
\begin{picture}(40,40)(-5,-8)  %graph for III(c_1)
    \thicklines
    \put(20,20){\circle*{3}}
    \put(0,20){\circle*{3}}
    \put(0,0){\circle*{3}}
    \put(20,0){\circle*{3}}
    \put(0,20){\line(1,0){20}}
    \put(0,0){\line(0,1){20}}
    \put(20,0){\line(-1,1){20}}
    \put(0,0){\line(1,1){20}}
    \put(20,0){\line(0,1){20}}
    \put(0,0){\line(1,0){20}}
    %\put(12,8){\vector(-1,1){0}}
    \put(20,16){\vector(0,1){0}}
    \put(4,4){\vector(-1,-1){0}}
    \put(15,20){\vector(1,0){0}}
    \put(0,16){\vector(0,1){0}}
    \put(16,4){\vector(1,-1){0}}
    \put(5,0){\vector(-1,0){0}}
  \end{picture} &
  \begin{picture}(40,40)(-5,-8)  %graph for III(c_1')
    \thicklines
    \put(20,20){\circle*{3}}
    \put(0,20){\circle*{3}}
    \put(0,0){\circle*{3}}
    \put(20,0){\circle*{3}}
    \put(0,20){\line(1,0){20}}
    \put(0,0){\line(0,1){20}}
    \put(20,0){\line(-1,1){20}}
    \put(0,0){\line(1,1){20}}
    \put(20,0){\line(0,1){20}}
    \put(0,0){\line(1,0){20}}
    %\put(12,8){\vector(-1,1){0}}
    \put(20,8){\vector(0,-1){0}}
    \put(4,4){\vector(-1,-1){0}}
    \put(15,20){\vector(1,0){0}}
    \put(0,16){\vector(0,1){0}}
    \put(16,4){\vector(1,-1){0}}
    \put(5,0){\vector(-1,0){0}}
  \end{picture} &
  \begin{picture}(40,40)(-5,-8)  %graph for III(c_1'')
    \thicklines
    \put(20,20){\circle*{3}}
    \put(0,20){\circle*{3}}
    \put(0,0){\circle*{3}}
    \put(20,0){\circle*{3}}
    \put(0,20){\line(1,0){20}}
    \put(0,0){\line(0,1){20}}
    \put(20,0){\line(-1,1){20}}
    \put(0,0){\line(1,1){20}}
    \put(20,0){\line(0,1){20}}
    \put(0,0){\line(1,0){20}}
    %\put(12,8){\vector(-1,1){0}}
    \put(20,8){\vector(0,-1){0}}
    \put(9,9){\vector(1,1){0}}
    \put(15,20){\vector(1,0){0}}
     \put(0,16){\vector(0,1){0}}
    \put(16,4){\vector(1,-1){0}}
    \put(5,0){\vector(-1,0){0}}
  \end{picture}\\
 \end{tabular}

\bigskip
\bigskip

Finally, for Type IV diagrams,   first of all, we may without loss of generality assume that $(d_2/p_3)=-1$. There are two families (I) and (II) of graphs depending on whether or not the subgraph on the points $d_1,d_2,d_3$ forms a cycle. Thus

\noindent (I) are those where the arrows on $\{d_1,d_2,d_3\}$  form (the unique) cycle:

$$
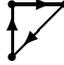
\begin{picture}(40,40)(-5,-8)
    \thicklines
    \put(20,20){\circle*{3}}
    \put(0,20){\circle*{3}}
    \put(0,0){\circle*{3}}
    %\put(20,0){\circle*{3}}
    \put(0,20){\line(1,0){20}}
    \put(0,0){\line(0,1){20}}
    %\put(20,0){\line(-1,1){20}}
    \put(0,0){\line(1,1){20}}
    %\put(20,0){\line(0,1){20}}
    %\put(0,0){\line(1,0){20}}
    %\put(12,8){\vector(-1,1){0}}
    %\put(20,16){\vector(0,1){0}}
    \put(4,4){\vector(-1,-1){0}}
    \put(15,20){\vector(1,0){0}}
    \put(0,16){\vector(0,1){0}}
    %\put(11,9){\vector(-1,1){0}}
    %\put(5,0){\vector(-1,0){0}}
    \end{picture} $$

\noindent (II) are those where the arrows do not form a cycle, from which we pick

$$
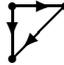
\begin{picture}(40,40)(-5,-8)
    \thicklines
     \put(20,20){\circle*{3}}
    \put(0,20){\circle*{3}}
    \put(0,0){\circle*{3}}
   %\put(20,0){\circle*{3}}
    \put(0,20){\line(1,0){20}}
    \put(0,0){\line(0,1){20}}
    %\put(20,0){\line(-1,1){20}}
    \put(0,0){\line(1,1){20}}
    %\put(20,0){\line(0,1){20}}
    %\put(0,0){\line(1,0){20}}
    %\put(12,8){\vector(-1,1){0}}
    %\put(20,16){\vector(0,1){0}}
    \put(4,4){\vector(-1,-1){0}}
    \put(15,20){\vector(1,0){0}}
    \put(0,8){\vector(0,-1){0}}
    %\put(11,9){\vector(-1,1){0}}
    %\put(5,0){\vector(-1,0){0}}
  \end{picture}$$

With respect to (I) there are four possible inequivalent diagrams depending on the number of line segments connecting this triangle to $d_4$. (Equivalence is determined by bijections of $\{d_1,d_2,d_3\}$ to itself.)
On the other hand there are eight diagrams (all inequivalent) corresponding to (II).  See Appendix I for a complete list of inequivalent diagrams.

\bigskip

\noindent{\em Remark.} The fields $k$ with discriminants of Types I and II for which $\Cl_2(k)$ is elementary may be found in \cite{BS}, but are given here in a much more concise package. We have also rearranged the sequences of discriminants ($a_j$) and ($b_j$), for $j=1,\dots,13$, so that first, they match up pairwise, i.e., $(a_j,b_j$), as closely  as possible, and moreover, so that the individual sequences, i.e., $\{(a_j)\}$ and  $\{(b_j)\}$,  are grouped with some of the  consecutive terms having similar relevant class group structures, as seen in Appendix II.

\newpage

\section{{\bf Appendix I\\$d_k=d_1d_2d_3d_4$}}

See Section~\ref{S2} for a description of the concepts and notation.
\bigskip

\noindent{\bf Type I}\;\;  ($d_1,d_2>0$, $d_3,d_4<0$, no $d_j=-4$)

\mbox{}
\medskip

\noindent $4$-$\rank=2$:

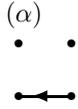
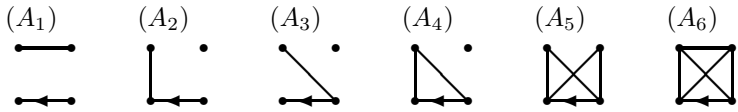
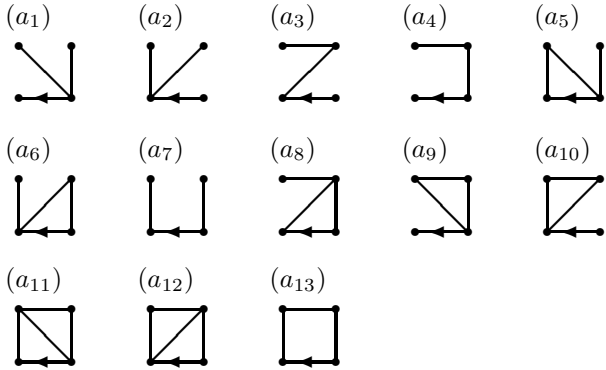
\mbox{}
\medskip

\noindent
\begin{tabular}{ll}\vspace{-.1in}
($\alpha$)\\
 \begin{picture}(40,40)(-5,-8)  %graph for I\alpha
    \thicklines
    \put(20,20){\circle*{3}}
    \put(0,20){\circle*{3}}
    \put(0,0){\circle*{3}}
    \put(20,0){\circle*{3}}
    % \put(0,20){\line(1,0){20}}
    \put(0,0){\line(1,0){20}}
    %\put(20,0){\line(-1,1){20}}
     \put(5,0){\vector(-1,0){0}}
  \end{picture}\\
\end{tabular}

\mbox
\medskip

\noindent $4$-$\rank=1$:

\mbox{}
\medskip

\noindent \begin{tabular}{llllll}\vspace{-.1in}
($A_1$)& ($A_2$) & ($A_3$) & ($A_4$) & ($A_5$) & ($A_6$)\\
\end{tabular}

\mbox{}
\medskip

\noindent $4$-$\rank=0$:

\mbox
\medskip

\noindent \begin{tabular}{lllll}\vspace{-.1in}
($a_1$)   & ($a_2$)   & ($a_3$)  & ($a_4$) & ($a_5$)  \\
\end{tabular}

\mbox
\medskip

\noindent \begin{tabular}{lllll}\vspace{-.1in}
($a_6$)   & ($a_7$)   & ($a_8$)  & ($a_9$) & ($a_{10}$)  \\
\end{tabular}
\mbox
\medskip

\noindent \begin{tabular}{lll}\vspace{-.1in}
($a_{11}$)   & ($a_{12}$)   & ($a_{13}$)   \\
\begin{picture}(40,40)(-5,-8)  %graph for I(a_11)
    \thicklines
    \put(20,20){\circle*{3}}
    \put(0,20){\circle*{3}}
    \put(0,0){\circle*{3}}
    \put(20,0){\circle*{3}}
    \put(0,20){\line(1,0){20}}
    \put(0,0){\line(1,0){20}}
    \put(0,0){\line(0,1){20}}
    \put(20,0){\line(-1,1){20}}
     %\put(0,0){\line(1,1){20}}
     \put(20,0){\line(0,1){20}}
     \put(5,0){\vector(-1,0){0}}
  \end{picture}&
\begin{picture}(40,40)(-5,-8)  %graph for I(a_12)
    \thicklines
    \put(20,20){\circle*{3}}
    \put(0,20){\circle*{3}}
    \put(0,0){\circle*{3}}
    \put(20,0){\circle*{3}}
    \put(0,20){\line(1,0){20}}
    \put(0,0){\line(1,0){20}}
    \put(0,0){\line(0,1){20}}
    %\put(20,0){\line(-1,1){20}}
     \put(0,0){\line(1,1){20}}
     \put(20,0){\line(0,1){20}}
     \put(5,0){\vector(-1,0){0}}
  \end{picture} &
\begin{picture}(40,40)(-5,-8)  %graph for I(a_13)
    \thicklines
    \put(20,20){\circle*{3}}
    \put(0,20){\circle*{3}}
    \put(0,0){\circle*{3}}
    \put(20,0){\circle*{3}}
    \put(0,20){\line(1,0){20}}
    \put(0,0){\line(1,0){20}}
    \put(0,0){\line(0,1){20}}
    %\put(20,0){\line(-1,1){20}}
    % \put(0,0){\line(1,1){20}}
     \put(20,0){\line(0,1){20}}
     \put(5,0){\vector(-1,0){0}}
  \end{picture}\\
\end{tabular}

\newpage

\noindent{\bf Type II}\;\;  ($d_1,d_2>0$, \;$d_3,d_4<0$, \;$d_4=-4$)

\mbox{}
\medskip

\noindent $4$-$\rank=2$:

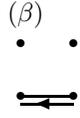
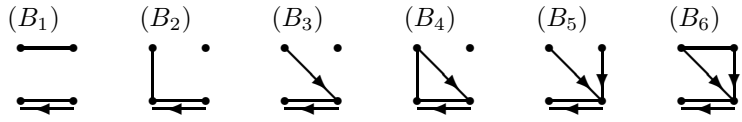
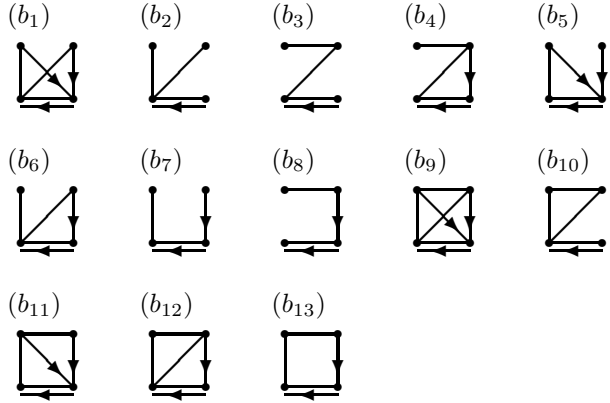
\mbox{}
\medskip

\noindent \begin{tabular}{l}\vspace{-.1in}
 ($\beta$)\\
 \begin{picture}(40,40)(-5,-8)  %graph for II(\beta)
    \thicklines
    \put(20,20){\circle*{3}}
    \put(0,20){\circle*{3}}
    \put(0,0){\circle*{3}}
    \put(20,0){\circle*{3}}
   % \put(0,20){\line(1,0){20}}
    \put(0,0){\line(1,0){20}}
    \put(0,-3){\line(1,0){20}}
    %\put(0,0){\line(0,1){20}}
    %\put(20,0){\line(-1,1){20}}
     %\put(0,0){\line(1,1){20}}
     %\put(20,0){\line(0,1){20}}
     \put(5,-3){\vector(-1,0){0}}
  \end{picture}\\
 \end{tabular}

\mbox{}
\medskip

\noindent $4$-$\rank=1$:

\mbox{}
\medskip

\noindent \begin{tabular}{llllll}\vspace{-.1in}
 ($B_1$) & ($B_2$) & ($B_3$) & ($B_4$) & ($B_5$) & ($B_6$)  \\
\end{tabular}

\mbox{}
\medskip

\noindent $4$-$\rank=0$:

\mbox{}
\medskip

\noindent \begin{tabular}{lllll}\vspace{-.1in}
($b_1$) & ($b_2$) & ($b_3$) & ($b_4$) & ($b_5$) \\
\end{tabular}

\mbox{}
\medskip

\noindent \begin{tabular}{lllll}\vspace{-.1in}
($b_6$) & ($b_7$) & ($b_8$) & ($b_9$) & ($b_{10}$) \\
\end{tabular}

\mbox{}
\medskip

\noindent \begin{tabular}{lll}\vspace{-.1in}
($b_{11}$) & ($b_{12}$) & ($b_{13}$) \\
\begin{picture}(40,40)(-5,-8)  %graph for II(b_11)
    \thicklines
    \put(20,20){\circle*{3}}
    \put(0,20){\circle*{3}}
    \put(0,0){\circle*{3}}
    \put(20,0){\circle*{3}}
    \put(0,20){\line(1,0){20}}
    \put(0,0){\line(1,0){20}}
    \put(0,-3){\line(1,0){20}}
     \put(0,0){\line(0,1){20}}
    \put(20,0){\line(-1,1){20}}
    % \put(0,0){\line(1,1){20}}
     \put(20,0){\line(0,1){20}}
     \put(5,-3){\vector(-1,0){0}}
      \put(16,4){\vector(1,-1){0}}
       \put(20,4){\vector(0,-1){0}}
  \end{picture}&
  \begin{picture}(40,40)(-5,-8)  %graph for II(b_12)
    \thicklines
    \put(20,20){\circle*{3}}
    \put(0,20){\circle*{3}}
    \put(0,0){\circle*{3}}
    \put(20,0){\circle*{3}}
   \put(0,20){\line(1,0){20}}
    \put(0,0){\line(1,0){20}}
    \put(0,-3){\line(1,0){20}}
   \put(0,0){\line(0,1){20}}
   % \put(20,0){\line(-1,1){20}}
     \put(0,0){\line(1,1){20}}
     \put(20,0){\line(0,1){20}}
     \put(5,-3){\vector(-1,0){0}}
      %\put(16,4){\vector(1,-1){0}}
       \put(20,4){\vector(0,-1){0}}
  \end{picture}&
  \begin{picture}(40,40)(-5,-8)  %graph for II(b_13)
    \thicklines
    \put(20,20){\circle*{3}}
    \put(0,20){\circle*{3}}
    \put(0,0){\circle*{3}}
    \put(20,0){\circle*{3}}
    \put(0,20){\line(1,0){20}}
    \put(0,0){\line(1,0){20}}
    \put(0,-3){\line(1,0){20}}
    \put(0,0){\line(0,1){20}}
    %\put(20,0){\line(-1,1){20}}
    % \put(0,0){\line(1,1){20}}
     \put(20,0){\line(0,1){20}}
     \put(5,-3){\vector(-1,0){0}}
      %\put(16,4){\vector(1,-1){0}}
       \put(20,4){\vector(0,-1){0}}
  \end{picture}\\
\end{tabular}

\newpage

\noindent{\bf Type III}\;\;  ($d_1,d_2,d_3,d_4<0$, \;$d_j\not=-4$ )

\mbox{}
\medskip

\noindent $4$-$\rank=2$:

\mbox{}
\medskip

None

\mbox{}
\medskip

\noindent $4$-$\rank=1$:

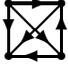
\mbox{}
\medskip

\noindent \begin{tabular}{l}\vspace{-.1in}
($C$) \\
\begin{picture}(40,40)(-5,-8)  %graph for III(C)
    \thicklines
    \put(20,20){\circle*{3}}
    \put(0,20){\circle*{3}}
    \put(0,0){\circle*{3}}
    \put(20,0){\circle*{3}}
    \put(0,20){\line(1,0){20}}
    \put(0,0){\line(0,1){20}}
    \put(20,0){\line(-1,1){20}}
    \put(0,0){\line(1,1){20}}
    \put(20,0){\line(0,1){20}}
    \put(0,0){\line(1,0){20}}
    %\put(12,8){\vector(-1,1){0}}
    \put(20,16){\vector(0,1){0}}
    \put(4,4){\vector(-1,-1){0}}
    \put(15,20){\vector(1,0){0}}
    \put(0,8){\vector(0,-1){0}}
    \put(16,4){\vector(1,-1){0}}
    \put(5,0){\vector(-1,0){0}}
  \end{picture}\\
 \end{tabular}

\mbox{}
\medskip

\noindent $4$-$\rank=0$:

\mbox{}
\medskip

\noindent \begin{tabular}{lll}\vspace{-.1in}
($c_1$)  &   ($c_2$)  &   ($c_3$)\\
\begin{picture}(40,40)(-5,-8)  %graph for III(c_1)
    \thicklines
    \put(20,20){\circle*{3}}
    \put(0,20){\circle*{3}}
    \put(0,0){\circle*{3}}
    \put(20,0){\circle*{3}}
    \put(0,20){\line(1,0){20}}
    \put(0,0){\line(0,1){20}}
    \put(20,0){\line(-1,1){20}}
    \put(0,0){\line(1,1){20}}
    \put(20,0){\line(0,1){20}}
    \put(0,0){\line(1,0){20}}
    %\put(12,8){\vector(-1,1){0}}
    \put(20,16){\vector(0,1){0}}
    \put(4,4){\vector(-1,-1){0}}
    \put(15,20){\vector(1,0){0}}
    \put(0,16){\vector(0,1){0}}
    \put(16,4){\vector(1,-1){0}}
    \put(5,0){\vector(-1,0){0}}
  \end{picture} &
  \begin{picture}(40,40)(-5,-8)  %graph for III(c_2)
    \thicklines
    \put(20,20){\circle*{3}}
    \put(0,20){\circle*{3}}
    \put(0,0){\circle*{3}}
    \put(20,0){\circle*{3}}
    \put(0,20){\line(1,0){20}}
    \put(0,0){\line(0,1){20}}
    \put(20,0){\line(-1,1){20}}
    \put(0,0){\line(1,1){20}}
    \put(20,0){\line(0,1){20}}
    \put(0,0){\line(1,0){20}}
    %\put(12,8){\vector(-1,1){0}}
    \put(20,16){\vector(0,1){0}}
    \put(9,9){\vector(1,1){0}}
    \put(15,20){\vector(1,0){0}}
    \put(0,16){\vector(0,1){0}}
    \put(16,4){\vector(1,-1){0}}
    \put(5,0){\vector(-1,0){0}}
  \end{picture} &
  \begin{picture}(40,40)(-5,-8)  %graph for III(c_3)
    \thicklines
    \put(20,20){\circle*{3}}
    \put(0,20){\circle*{3}}
    \put(0,0){\circle*{3}}
    \put(20,0){\circle*{3}}
    \put(0,20){\line(1,0){20}}
    \put(0,0){\line(0,1){20}}
    \put(20,0){\line(-1,1){20}}
    \put(0,0){\line(1,1){20}}
    \put(20,0){\line(0,1){20}}
    \put(0,0){\line(1,0){20}}
    %\put(12,8){\vector(-1,1){0}}
    \put(20,8){\vector(0,-1){0}}
    \put(9,9){\vector(1,1){0}}
    \put(15,20){\vector(1,0){0}}
    \put(0,8){\vector(0,-1){0}}
    \put(16,4){\vector(1,-1){0}}
    \put(5,0){\vector(-1,0){0}}
  \end{picture}\\
 \end{tabular}

 \bigskip
 \bigskip
 \bigskip
 \bigskip

 \noindent{\bf Type IV}\;\;  ($d_1,d_2,d_3,d_4<0$, \;$d_4=-4$ )

\mbox{}
\medskip

\noindent $4$-$\rank=2$:

\mbox{}
\medskip

None

\mbox{}
\medskip

\noindent $4$-$\rank=1$:

\mbox{}
\medskip

\noindent \begin{tabular}{llll}\vspace{-.1in}
($D_1$)  &   ($D_2$)  &   ($D_3$) & ($D_4$) \\
\end{tabular}

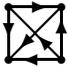
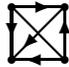
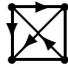
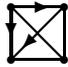
 \mbox{}
\medskip

\noindent $4$-$\rank=0$:

\mbox{}
\medskip

\noindent \begin{tabular}{llll}\vspace{-.1in}
($d_1$)  &   ($d_2$)  &   ($d_3$) & ($d_4$) \\
\end{tabular}
\mbox{}
\medskip

\noindent \begin{tabular}{llll}\vspace{-.1in}
($d_5$)  &   ($d_6$)  &   ($d_7$) & ($d_8$) \\
\end{tabular}

\newpage

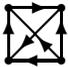
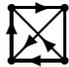
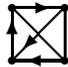
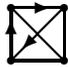
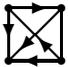
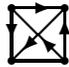
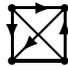
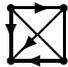
\begin{center} {\bf Appendix I$'$} \end{center}

\bigskip

Here we give data equivalent to that in Appendix~I by listing the values of the relevant Kronecker symbols for each case.

\bigskip

\noindent{\bf Type I}\;\;  ($d_1,d_2>0$, $d_3,d_4<0$, no $d_j=-4$)

\medskip
For each case, $(d_3/p_4)=+1$ and

\medskip

\begin{tabular}{c|ccccc}
   \text{case} & $(d_1/p_2)$ & $(d_1/p_3)$ & $(d_2/p_3)$
            & $(d_1/p_4)$ & $(d_2/p_4)$ \\ \hline
   \rsp $(\alpha)$ & $+1$ & $+1$ & $+1$ & $+1$ & $+1$ \\
   \rsp $(A_1)$ & $-1$ & $+1$ & $+1$ & $+1$ & $+1$  \\
   \rsp $(A_2)$ & $+1$ & $-1$ & $+1$ & $+1$ & $+1$   \\
   \rsp $(A_3)$ & $+1$ & $+1$ & $+1$ & $-1$ & $+1$  \\
   \rsp $(A_4)$ & $+1$ & $-1$ & $+1$ & $-1$ & $+1$ \\
   \rsp $(A_5)$ & $+1$ & $-1$ & $-1$ & $-1$ & $-1$  \\
   \rsp $(A_6)$ & $-1$ & $-1$ & $-1$ & $-1$ & $-1$ \\
   \rsp $(a_1)$ & $+1$ & $+1$ & $+1$ & $-1$ & $-1$ \\
   \rsp $(a_2)$ & $+1$ & $-1$ & $-1$ & $+1$ & $+1$  \\
   \rsp $(a_3)$ & $-1$ & $+1$ & $-1$ & $+1$ & $+1$   \\
   \rsp $(a_4)$ & $-1$ & $+1$ & $+1$ & $+1$ & $-1$  \\
   \rsp $(a_5)$ & $+1$ & $-1$ & $+1$ & $-1$ & $-1$ \\
   \rsp $(a_6)$ & $+1$ & $-1$ & $-1$ & $+1$ & $-1$  \\
   \rsp $(a_7)$ & $+1$ & $-1$ & $+1$ & $+1$ & $-1$ \\
   \rsp $(a_8)$ & $-1$ & $+1$ & $-1$ & $+1$ & $-1$  \\
   \rsp $(a_9)$ & $-1$ & $+1$ & $+1$ & $-1$ & $-1$   \\
   \rsp $(a_{10})$ & $-1$ & $-1$ & $-1$ & $+1$ & $+1$  \\
   \rsp $(a_{11})$ & $-1$ & $-1$ & $+1$ & $-1$ & $-1$ \\
   \rsp $(a_{12})$ & $-1$ & $-1$ & $-1$ & $+1$ & $-1$  \\
   \rsp $(a_{13})$ & $-1$ & $-1$ & $+1$ & $+1$ & $-1$ \\
  \end{tabular}

\newpage

\noindent{\bf Type II}\;\;  ($d_1,d_2>0$, $d_3,d_4<0$,  $d_4=-4$)

\medskip
For each case, $(d_3/2)=\pm 1$ and

\medskip

\begin{tabular}{c|ccccc}
   \text{case} & $(d_1/p_2)$ & $(d_1/p_3)$ & $(d_2/p_3)$
            & $(d_1/2)$ & $(d_2/2)$ \\ \hline
   \rsp $(\beta)$ & $+1$ & $+1$ & $+1$ & $+1$ & $+1$ \\
   \rsp $(B_1)$ & $-1$ & $+1$ & $+1$ & $+1$ & $+1$  \\
   \rsp $(B_2)$ & $+1$ & $-1$ & $+1$ & $+1$ & $+1$   \\
   \rsp $(B_3)$ & $+1$ & $+1$ & $+1$ & $-1$ & $+1$  \\
   \rsp $(B_4)$ & $+1$ & $-1$ & $+1$ & $-1$ & $+1$ \\
   \rsp $(B_5)$ & $+1$ & $+1$ & $+1$ & $-1$ & $-1$  \\
   \rsp $(B_6)$ & $-1$ & $+1$ & $+1$ & $-1$ & $-1$ \\
   \rsp $(b_1)$ & $+1$ & $-1$ & $-1$ & $-1$ & $-1$ \\
   \rsp $(b_2)$ & $+1$ & $-1$ & $-1$ & $+1$ & $+1$  \\
   \rsp $(b_3)$ & $-1$ & $+1$ & $-1$ & $+1$ & $+1$   \\
   \rsp $(b_4)$ & $-1$ & $+1$ & $-1$ & $+1$ & $-1$  \\
   \rsp $(b_5)$ & $+1$ & $-1$ & $+1$ & $-1$ & $-1$ \\
   \rsp $(b_6)$ & $+1$ & $-1$ & $-1$ & $+1$ & $-1$  \\
   \rsp $(b_7)$ & $+1$ & $-1$ & $+1$ & $+1$ & $-1$ \\
   \rsp $(b_8)$ & $-1$ & $+1$ & $+1$ & $+1$ & $-1$  \\
   \rsp $(b_9)$ & $-1$ & $-1$ & $-1$ & $-1$ & $-1$   \\
   \rsp $(b_{10})$ & $-1$ & $-1$ & $-1$ & $+1$ & $+1$  \\
   \rsp $(b_{11})$ & $-1$ & $-1$ & $+1$ & $-1$ & $-1$ \\
   \rsp $(b_{12})$ & $-1$ & $-1$ & $-1$ & $+1$ & $-1$  \\
   \rsp $(b_{13})$ & $-1$ & $-1$ & $+1$ & $+1$ & $-1$ \\
  \end{tabular}

\newpage

\noindent{\bf Type III}\;\;  ($d_1,d_2,d_3,d_4<0$, \;$d_j\not=-4$ )

\medskip
$(d_1/p_2)=(d_1/p_4)=(d_4/p_3)=-1$ and

\medskip

\begin{tabular}{c|ccc}
   \text{case} & $(d_1/p_3)$ & $(d_2/p_3)$ & $(d_2/p_4)$  \\ \hline
   \rsp $(C)$ & $-1$ & $-1$ & $+1$ \\
   \rsp $(c_1)$ & $+1$ & $-1$ & $+1$  \\
   \rsp $(c_2)$ & $+1$ & $+1$ & $+1$   \\
   \rsp $(c_3)$ & $-1$ & $+1$ & $-1$   \\
\end{tabular}

 \bigskip
 \bigskip

 \noindent{\bf Type IV}\;\;  ($d_1,d_2,d_3,d_4<0$, \;$d_4=-4$ )

\medskip

$(d_1/p_2)=(d_2/p_3)=-1$

\medskip
\begin{tabular}{c|cccc}
   \text{case} & $(d_1/p_3)$ & $(d_1/2)$ & $(d_2/2)$
            & $(d_3/2)$  \\ \hline
   \rsp $(D_1)$ & $-1$ & $+1$ & $+1$ & $+1$   \\
   \rsp $(D_2)$ & $-1$ & $-1$ & $+1$ & $+1$   \\
   \rsp $(D_3)$ & $-1$ & $+1$ & $-1$ & $-1$   \\
   \rsp $(D_4)$ & $-1$ & $-1$ & $-1$ & $-1$   \\
   \rsp $(d_1)$ & $+1$ & $+1$ & $+1$ & $+1$  \\
   \rsp $(d_2)$ & $+1$ & $+1$ & $-1$ & $+1$  \\
   \rsp $(d_3)$ & $+1$ & $-1$ & $-1$ & $+1$  \\
   \rsp $(d_4)$ & $+1$ & $-1$ & $-1$ & $-1$  \\
   \rsp $(d_5)$ & $-1$ & $+1$ & $-1$ & $+1$  \\
   \rsp $(d_6)$ & $-1$ & $+1$ & $+1$ & $-1$   \\
   \rsp $(d_7)$ & $-1$ & $-1$ & $+1$ & $-1$  \\
   \rsp $(d_8)$ & $-1$ & $-1$ & $-1$ & $+1$   \\
  \end{tabular}

\newpage

\section{{\bf Review of the Structure of ${\Gal(k^2/k)}$ when ${\Cl_2(k)\simeq (2,2)}$}}\label{S4}

We  now review the structure of $G=\Gal(k^2/k)$ for  discriminants $d_k$ of Types I and II given above and for which $\Cl_2(k)\simeq (2,2)$. When all four prime discriminants are negative, i.e., for $d_k$ of Types III and IV (including the non-elementary cases) we know that $G$ is abelian, by Theorem~1 of \cite{BS} or Theorem~2 in \cite{BLS}. Thus $G\simeq \Cl_2(k)$ for these two types.

Recall that if $G$ is a finite 2-group whose abelianization $G^{\ab}=G/G'\simeq (2,2)$, where $G'$ is the commutator subgroup of $G$,  then $G$ is isomorphic to a member of  one of the following families of groups:
\bigskip

$V_4=(2,2)$, i.e., the Klein 4-group;

$Q_m$, the quaternion group of order $m\geq 8$;

$D_m$, the dihedral group of order $2m\geq 8$;

$SD_m$, the semi-dihedral group of order $m\geq 16$.
\medskip

We will also use the following abbreviations for convenience:
$$Q=Q_8\,,\;\;Q_g=Q_m\,(m\geq 16)\,,\;\;D=D_m\,,\;\;S=SD_m.$$

For these groups, let $H_1, \,H_2,\, H_3$ be the three maximal subgroups of $G$, where $H_1$ is cyclic in all the above families. (This subgroup is unique if $G$ is either $Q_g,\,D,$ or $S$. When $G$ is $(2,2)$ or $Q$, then all three are cyclic.)

Then recall that when $\Cl_2(k)\simeq (2,2)$,  we can determine which of the above groups $G=\Gal(k^2/k)$ is isomorphic to by computing the orders of the capitulation kernels in the three unramified quadratic extensions of $k$ along with Taussky's conditions $A$ and $B$ (see, for example, \cite{Kis}):

If $L/K$ is a finite extension of number fields, then define the transfer map by $j_{L/K}:\Cl_2(K)\rightarrow \Cl_2(L)$ as the homomorphism of classes induced by extension of ideals from $K$ to $L$. Then the capitulation kernel of $L/K$ is the kernel $\ker j_{L/K}$. Now suppose $L/K$ is an unramified quadratic extension. Then the Taussky conditions are given by

\medskip
$L/K$ satisfies condition $A$, if $|\ker j_{L/K}\cap N_{L/K} \Cl(L)|>1,$

$L/K$ satisfies condition $B$, if $|\ker j_{L/K}\cap N_{L/K} \Cl(L)|=1.$

\bigskip

In our case, we then have the following data: Let $K_i=(k^1)^{H_i}$ be the quadratic extension of $k$ in $k^1$ fixed by the maximal  subgroup $H_i$ of $G=\Gal(k^2/k)$ and denote $j_{K_i/k}$ by $j_i$ for $i=1,2,3$ (where again $H_1$ is cyclic), cf.\,\cite{Kis}.

\medskip
\begin{tabular}{c|cc}
$G$ & $|\ker j_1|\,(A/B)$ & $|\ker j_i|\,(A/B), (i=2,3)$\\\hline
\rsp $(2,2)$ & $4$  & $4$ \\
\rsp $Q$ & $2A$ & $2A$ \\
\rsp $Q_g$ & $2A$ & $2B$ \\
\rsp $D$ & $4$ & $2B$\\
\rsp $S$ & $2B$ & $2B$ \\
\end{tabular}

 In Appendix~II we list all the groups $G=\Gal(k^2/k)$ including their orders for each of the cases $(a_i)$ and $(b_i)$ given in Appendix~I. This is a much more concise version of some of the data given in \cite{BS}. We have also corrected a couple of errors occurring  in the tables of \cite{BS}.

\noindent{\em Remark.} We also note a couple of errors in Theorem~3 in \cite{BS} in the ambiguous case of $G=D$ or $Q_g$. On p.\,173. in part (2), when $G=D$, $K_c=k(\sqrt{ab}\,)$, not $k(\sqrt{a}\,)$. Also on p.\,174, in part (3), when $G=Q_g$, $K_c=k(\sqrt{a}\,)$, not $k(\sqrt{b}\,)$.

\bigskip

 Here's a description of the tables in Appendix~II.

Column 1: This column indicates the cases $(a_i)$ and $(b_i)$ as given in Appendix~I.

For the moment, we skip to Column 4.

Column 4: Shows the possible values of $N\eps_{12}$, the norm of the fundamental unit $\eps_{12}$ in $\Q(\sqrt{d_1d_2}\,)$.
When there are two possible values in any entry, we represent them in a $2$-d column vector given in rounded brackets. If there are several such column vectors in a given case, all their rows correspond to the same relevant data.

Column 2: Row 1 gives the relevant values $\nu_{ij}\in\{0,1\}$ such that $(-1)^{\nu_{ij}}=(d_i/p_j)$, where $p_i$ is the rational prime dividing $d_i$. Moreover, in the cases $(b_i)$ we include the two possible values of $\nu_{34}$, which is given as a column vector in square brackets. (Observe here that the only changes occurring are the values of $\delta$, $\delta_1$, and $\delta_2$ in Column~3.) Row 2 gives $N_i=N_{k_i/k}\Cl_2(k_i)$ for $k_1=k(\sqrt{d_1}\,)$, $k_2=k(\sqrt{d_2}\,)$, and $k_3=k(\sqrt{d_1d_2}\,)$.  Row 3 gives the unit indices  $q(k_i)=q(k_i/\Q)$, which is the index of $E_{k_i}$ modulo the subgroup generated by the unit groups in the quadratic subfields of $k_i$. (Computing these indices requires  data in the first row of the third column.)

Column 3: Row 1 gives the values $\delta=\delta(\varepsilon_k)=\sfk N(1+\varepsilon_k)$, where $\sfk$ denotes the square-free kernel. Similarly, $\delta_1=\delta(\eps_{234})$ and $\delta_2=\delta(\eps_{134})$. Here $\eps_{i34}$ denotes the fundamental unit of $\Q(\sqrt{d_id_3d_4}\,)$. We also observe that $\delta_3=\delta(\eps_{34})=p_4$ in the $(a_i)$ cases and in the $(b_i)$ cases $\delta(\eps_{34})=2p_3$ or $2$, according as $p_3\equiv 3$ or $7\bmod 8$, respectively, cf.\,\cite{Kub}.  Rows 2 and 3 are self-explanatory.
(Here $h_2(k_i)$ denotes the $2$-class number of $k_i$, etc. Also we write $h_2(\Q(\sqrt{m}\,))$ as $h_2(m).$)

Column 5: Is self-explanatory.

Column 6: Gives the quotient group $G^+/G_3^+$ for  $G^+=\Gal(k_+^2/k)$  for the narrow second Hilbert $2$-class field of $k$. This will be discussed in Section~7 below.

The methods for computing the data given in Columns 1-5 can be found in \cite{BS}. Moreover, $h_2(k_i)$ can be calculated using Kuroda's class number formula: $h_2(k_i)=\frac1{4}q(k_i/\Q)\prod_F h_2(F)$ where $F$ ranges over the three quadratic subfields of $k_i$, cf.\,\cite{Lem2}.

Finally, for the $(b_i)$ cases, $\fq$ denotes the prime ideal in $k$ above $2$.

\bigskip

We now give an example of how the data on $G$ in Appendix II are computed, by computing the data in case $(a_8)$. Here are the steps we follow:

\medskip

1. Determine $\delta, \delta_1, \delta_2$, using the genus characters on $k,\Q(\sqrt{d/d_1}\,), \Q(\sqrt{d/d_2}\,)$.

2. Determine $\Cl_2(k)$.

3. Calculate $E_i=E_{k_i}$ and then $q(k_i)$ for $i=1,2,3.$

4. Determine $NE_i=N_{k_i/k}E_i$, and then $\ker j_i$.

5. Calculate $N_i=N_{k_i/k}\Cl_2(k_i)$.

6. Use steps 4 and 5 to determine the Type of $G$.

7. Determine $h_2(k_i)$ by Kuroda's class number formula.

8. Get the order of $G$.

\newpage
Consider $k$ in case $(a_8)$. Then $d=d_k=d_1d_2d_3d_4\not\equiv 4\bmod 8$  for prime discriminants $d_i$ with $d_1,d_2>0$; $d_3,d_4<0$ such that
$$ (d_1/p_2)=(d_2/p_3)=(d_2/p_4)=(d_4/p_3)=-1,\;(d_1/p_3)=(d_1/p_4)=+1,  $$
where $p_i$ is the prime dividing $d_i$.
\medskip

1. We determine $\delta=\delta(\eps_k)$, $\delta_1=\delta(\eps_{234})$, and $\delta_2=\delta(\eps_{134})$ (where $\varepsilon_m$ is the fundamental unit of the relevant quadratic field). We start with $\delta$ in $k$; let $\chi_i$ be  the genus character associated with the prime discriminant $d_i$ for $i=1,2,3,4$. Hence if $\fp_j$ is the (ramified) prime ideal in $k$ containing $p_j\mid d_j$, then $$\chi_i(\fp_j)=\left\{
                                                                           \begin{array}{ll}
                                                                             (d_i/p_j), & \hbox{if}\;\; j\not=i, \\
                                                                             \big((d_k/d_i)/p_i\big), & \hbox{otherwise.}
                                                                           \end{array}
                                                                         \right.$$
Then in case $(a_8)$ be have the following table of values for the four genus characters of $k$:

\begin{tabular}{c|cccc}
$\fp$ & $\chi_1$ & $\chi_2$ & $\chi_3$ & $\chi_4$ \\\hline
\rsp $\fp_1$ & $-1$  & $-1$ & $+1$ & $+1$ \\
\rsp $\fp_2$ & $-1$  & $-1$ & $-1$ & $-1$ \\
\rsp $\fp_3$ & $+1$  & $-1$ & $+1$ & $-1$\\
\rsp $\fp_4$ & $+1$  & $-1$ & $+1$ & $-1$ \\
\end{tabular}.
But notice that $\fp_3\fp_4$ is the only product of the $\fp_i$ that has value $=+1$ for each genus characters. Thus $\delta=N\fp_3\fp_4=p_3p_4.$  For all the details on genus characters, cf.\;Chapter 2 of \cite{Lem1}.

Now consider $\delta_1$. Then the procedure is similar, except we consider the genus characters of $F=\Q(\sqrt{d_2d_3d_4}\,)$. There are three such characters $\chi_i$ associated with the prime discriminants $d_i$ for $i=2,3,4$. Then the table of values for the three ramified primes $\fp_j$ in $F$ is given by

\begin{tabular}{c|ccc}
$\fp$ & $\chi_2$ & $\chi_3$ & $\chi_4$\\\hline
\rsp $\fp_2$ & $+1$  & $-1$ & $-1$ \\
\rsp $\fp_3$ & $-1$  & $+1$ & $-1$ \\
\rsp $\fp_4$ & $-1$  & $+1$ & $-1$ \\
\end{tabular}.
Therefore, $\delta_1=p_3p_4$. (Again refer to \cite{Lem1} for the details.)

Finally consider $\delta_2$. This time we consider $F=\Q(\sqrt{d_1d_3d_4}\,)$ with discriminant $d_F=d_1d_3d_4$. The three genus characters are $\chi_i$ associated with $d_i$ for $i=1,3,4$.  The table of values is:

\begin{tabular}{c|ccc}
$\fp$ & $\chi_1$ & $\chi_3$ & $\chi_4$\\\hline
\rsp $\fp_1$ & $+1$  & $+1$ & $+1$ \\
\rsp $\fp_3$ & $+1$  & $-1$ & $-1$ \\
\rsp $\fp_4$ & $+1$  & $+1$ & $+1$ \\
\end{tabular}. Hence there are several possibilities for the value of $\delta_2$, namely, $p_1, p_4, p_1p_4$. (We observe that $\Cl_2(F)$ is not elementary here, as opposed to the two other quadratic fields above.) This ambiguity will allow for a couple of possible types of $G$, as we'll see.

\medskip

2.  Now we compute $\Cl_2(k)$. We know that $\Cl_2(k)\simeq (2,2)$ and since $d_k$ is divisible by negative prime discriminants, $\Cl_2(k)$ is generated by those ideal classes containing a ramified prime. But first notice that $\fp_1\fp_2\fp_3\fp_4\sim (\sqrt{d_k}\,)\sim 1.$  Moreover, $\fp_3\fp_4\sim (1+\eps_k)\sim 1$. Thus $\fp_1\sim\fp_2$ and $\fp_3\sim\fp_4$ and so we have $\Cl_2(k)=\la [\fp_1],[\fp_3]\ra.$

\medskip

3. Now we determine the unit groups $E_i$ of the fields $k_i$ for $i=1,2,3$. Consider $E_1$ where $k_1=\Q(\sqrt{d_1},\sqrt{d_2d_3d_4}\,)$. Hence $E_1\supseteq \la -1,\eps_1,\eps_{234},\eps_k\ra$. But since $\delta \delta_1$ is a square in $k_1$, we know that $\sqrt{\eps_k\eps_{234}}\in E_1$, cf.\;\cite{Kub}. Thus $E_1=\la -1,\eps_1,\eps_{234},\sqrt{\eps_k\eps_{234}}\ra$ and so $q(k_1)=2$.  Now for the moment skip $i=2$ and consider $E_3$, the unit group of $k_3=\Q(\sqrt{d_1d_2},\sqrt{d_3d_4}\,)$. Hence $E_3\supseteq \la -1,\eps_{12},\eps_{34},\eps_k\ra.$ But since $\delta$ is a square in $k_3$, we have $E_3=\la -1,\eps_{12},\eps_{34}, \sqrt{\eps_k}\,\ra,$ and thus $q(k_3)=2$. Now for $i=2$. Then $E_2\supseteq \la -1,\eps_2,\eps_{134}, \eps_k\ra.$ What happens here depends on the value of $\delta_2$. If $\delta_2=p_1$, then $\delta\delta_2$ is a square in $k_2$; but if $\delta_2=p_4$ or $p_1p_4$, then no nontrivial combination of the $\delta$'s is a square, thus we have
$$E_2=\left\{
        \begin{array}{ll}
          \la -1,\eps_2,\eps_{134}, \sqrt{\eps_k\eps_{134}}\,\ra, & \hbox{if}\; \delta_2=p_1, \\
          \la -1,\eps_2,\eps_{134}, \eps_k\ra, & \hbox{if}\; \delta_2\in \{p_4,p_1p_4\}.
        \end{array}
      \right.$$
Moreover, $q(k_2)=2$ or $1$ according as $\delta_2=p_1$ or $\in \{p_4,p_1p_4\},$ respectively.

\medskip

4. First notice that $NE_1=E_k=NE_3$ and so $(E_k:NE_i)=1$ for $i=1,3$. Also, if $\delta_2=p_1$, then $NE_2=E_k$; whereas if $\delta_2=p_4$ or $p_1p_4$, then $NE_2=\la -1,\eps_k^2\ra$, in which case $(E_k:NE_2)=2$.  Now recall that $|\ker j_i|=2(E_k:NE_i)$. Thus $|\ker j_i|=2$ for $i=1,3$ and $i=2$ except when $\delta_2\not=p_1$, in which case $|\ker j_2|=4$. Now we can actually determine $\ker j_i$.  First, notice that $\fp_1\OO_{k_1}=(\sqrt{p_1}\,)$ and thus $\ker j_1=\la [\fp_1]\ra$.
Similarly, if $\delta_2=p_1$, then $\ker j_2=\la [\fp_2]\ra=\la [\fp_1]\ra$ (since $\fp_2\sim\fp_1$ from above). If $\delta_2\not=p_1$, then $\ker j_2=\Cl_2(k).$ Finally, $\fp_3\OO_{k_3}\sim 1$, since the ideal in $F=\Q(\sqrt{d_3d_4}\,)$ containing $p_3$ is principal, being of order $2$ in $F$, which has odd class number. Therefore, $\ker j_3=\la [\fp_3]\ra.$

\medskip

5. Consider $N_1=N_{k_1/k}\Cl_2(k_1)=\{c\in \Cl_2(k): (c, k_1/k)=+1\}$, i.e., the kernel of the Artin symbol $(\;\cdot\;,k_1/k).$  But we have $\big([\fp_3],k(\sqrt{d_1}\,)/k\big)=(d_1/p_3)=+1.$ Hence $[\fp_3]\in N_1$ and since the order of $N_1$ is $2$, $N_1=\la [\fp_3]\,\ra.$ Similarly, notice that $([\fp_1\fp_3],k_2/k)=(d_2/p_1p_3)=+1$ and thus $N_2=\la [\fp_1\fp_3]\ra.$ Furthermore, $([\fp_1,k_3/k)=(d_3d_4/p_1)=+1$, implying that $N_3=\la [\fp_1]\,\ra$.

\medskip

6.  Now we can determine the Type of $G$. If $\delta_2=p_1$, then the order of each $\ker j_i$ is $2$. Moreover, we can see, by comparing $N_i$ with $\ker j_i$ for $i=1,2,3$, that  all the Taussky conditions are $2B$. Hence $G$ is semi-dihedral.
On the other hand, when $\delta_2=p_4$ or $p_1p_4$, then $\ker j_2$ has order $4$ (while the rest have order $2$). Therefore, $G$ is dihedral.

\medskip

7. Now we calculate $h_2(k_i)$.  By Kuroda's class number formula and since $k_1=\Q(\sqrt{d_1},\sqrt{d_2d_3d_4}\,)$,   $h_2(k_1)=\frac1{4}q(k_1)h_2(k)h_2(d_1)h_2(d_2d_3d_4)=4$, since $q(k_1)=2,\, h_2(k)=4,\, h_2(d_1)=1$, and $h_2(d_2d_3d_4)=2$ (for instance by considering $C_4$-splittings of $d_2d_3d_4$). Similarly, $h_2(k_3)=4$.  Now notice that \\ $h_2(k_2)=\frac1{4}q(k_2)h_2(k)h_2(d_2)h_2(d_1d_3d_4)=q(k_2)h_2(d_1d_3d_4).$ Therefore,
$$h_2(k_2)=\left\{
             \begin{array}{ll}
               2h_2(d_1d_3d_4), & \hbox{if}\; \delta_2=p_1, \\
               h_2(d_1d_3d_4), & \hbox{otherwise.}
             \end{array}
           \right. $$

\medskip

8. Now $$|G|=\left\{
             \begin{array}{ll}
               4h_2(d_1d_3d_4), & \hbox{if}\; \delta_2=p_1, \\
               2h_2(d_1d_3d_4), & \hbox{otherwise,}
             \end{array}
           \right. $$
since here $|G|=2h_2(k_2)$.

\newpage
\section{{\bf Appendix II}}

See Sections~\ref{S4} and \ref{S7} for a description of the concepts and notation.
\bigskip

{\tiny $$\begin{array}{|c|c|c|c|c|c|}\hline
\rsp {\rm Case} &  (\nu_{12},\nu_{13},\nu_{23},\nu_{14},\nu_{24},\nu_{34})&\delta\qquad \delta_1\qquad \delta_2 & N\varepsilon_{12}& \Cl_2(k) &\\
  & N_1\;\;\;\;N_2\;\;\;\;N_3 & \ker j_1\;\;\;\ker j_2\;\;\;\ker j_3  &  & \text{Type of $G$} & G^+/G^+_3  \\
  & q(k_1)\;\;q(k_2)\;\; q(k_3) & h_2(k_1)\;\;\;h_2(k_2)\;\;\;h_2(k_3)&  &  \text{Order of $G$} &\\\hline\hline
a_1 &(0,0,0,1,1,0) & p_1p_2p_4\;\;\;\;\;\;p_2p_4\;\;\;\;\;\;p_1p_4 & \left(
                                                             \begin{array}{c}
                                                               -1 \\
                                                               +1 \\
                                                             \end{array}
                                                           \right)      & \la [\fp_1],[\fp_2]\ra & \\
& \la [\fp_2]\ra\;\; \la [\fp_1]\ra \;\;\la [\fp_1\fp_2]\ra & \la [\fp_1]\ra \;\;\; \la [\fp_2]\ra \;\;\; \left(
                                                                                                 \begin{array}{c}
                                                                                                   \la [\fp_1\fp_2]\ra \\
                                                                                                   \la [\fp_1],[\fp_2]\ra   \\
                                                                                                 \end{array}
                                                                                               \right) & & \left(
                                                                                                             \begin{array}{c}
                                                                                                               Q_g \\
                                                                                                               D \\
                                                                                                             \end{array}
                                                                                                           \right)
                                                                                        & 64.144        \\
 & 2\qquad 2\qquad 2 & 4\qquad 4\qquad 2h_2(d_1d_2) & & 4h_2(d_1d_2) &\\\hline
a_2 &(0,1,1,0,0,0) & p_4\;\;\;\;\;\;p_4\;\;\;\;\;\;p_4 & \left(
                                                             \begin{array}{c}
                                                               -1 \\
                                                               +1 \\
                                                             \end{array}
                                                           \right)      & \la [\fp_1],[\fp_2]\ra  & \\
& \la [\fp_2]\ra\;\; \la [\fp_1]\ra \;\;\la [\fp_1\fp_2]\ra & \la [\fp_1]\ra \;\;\; \la [\fp_2]\ra \;\;\; \left(
                                                                                                 \begin{array}{c}
                                                                                                   \la [\fp_1\fp_2]\ra \\
                                                                                                   \la [\fp_1],[\fp_2]\ra   \\
                                                                                                 \end{array}
                                                                                               \right) & & \left(
                                                                                                             \begin{array}{c}
                                                                                                               Q_g \\
                                                                                                               D \\
                                                                                                             \end{array}
                                                                                                           \right)
                                                                                           & 64.144     \\
 & 2\qquad 2\qquad 2 & 4\qquad 4\qquad 2h_2(d_1d_2) & & 4h_2(d_1d_2) &\\\hline
a_3 &(1,0,1,0,0,0) & p_4 \;\;\; p_4 \:\:\: \left(
                                                             \begin{array}{c}
                                                               p_4 \\
                                                               p_1 \;\text{or}\;p_1p_4 \\
                                                             \end{array}
                                                           \right)     & -1   & \la [\fp_1],[\fp_2]\ra & \\
& \la [\fp_1\fp_2]\ra\;\; \la [\fp_2]\ra \;\;\la [\fp_1]\ra & \la [\fp_1]\ra \;\;\; \left(
                                                                \begin{array}{c}
                                                                  \la [\fp_2]\ra \\
                                                                   \la [\fp_1],[\fp_2]\ra \\
                                                                \end{array}
                                                              \right)
 \;\;\; \la [\fp_1\fp_2]\ra  & & \left(
                                                                                                             \begin{array}{c}
                                                                                                               Q_g \\
                                                                                                               D \\
                                                                                                             \end{array}
                                                                                                           \right)
                                                                                       & 64.144         \\
 & 2\quad \left(
     \begin{array}{c}
       2 \\
       1 \\
     \end{array}
   \right)
   \quad 2 & 4\quad \left(
                       \begin{array}{c}
                         2h_2(d_1d_3d_4) \\
                         h_2(d_1d_3d_4) \\
                       \end{array}
                     \right)
   \quad 4  & & \left(
                       \begin{array}{c}
                         4h_2(d_1d_3d_4) \\
                         2h_2(d_1d_3d_4) \\
                       \end{array}
                     \right) & \\\hline
 a_4 &(1,0,0,0,1,0) & p_1p_2p_4\;\;\; p_2p_4 \;\;\;\left(
                                                             \begin{array}{c}
                                                               p_1p_4 \\
                                                               p_1 \;\text{or}\;p_4 \\
                                                             \end{array}
                                                           \right)     & -1   & \la [\fp_1],[\fp_2]\ra & \\
& \la [\fp_1\fp_2]\ra\;\; \la [\fp_2]\ra \;\;\la [\fp_1]\ra & [\fp_1] \;\;\;\left(
                                                                \begin{array}{c}
                                                                  \la [\fp_2]\ra \\
                                                                   \la [\fp_1],[\fp_2]\ra \\
                                                                \end{array}
                                                              \right)
 \;\;\; \la [\fp_1\fp_2]\ra   & & \left(
                                                                                                             \begin{array}{c}
                                                                                                               Q_g \\
                                                                                                               D \\
                                                                                                             \end{array}
                                                                                                           \right)
                                                                                           &64.144     \\
 & 2 \quad \left(
     \begin{array}{c}
       2 \\
       1 \\
     \end{array}
   \right)
   \quad 2 & 4\quad  \left(
                       \begin{array}{c}
                         2h_2(d_1d_3d_4) \\
                         h_2(d_1d_3d_4) \\
                       \end{array}
                     \right)
   \quad 4 & & \left(
                       \begin{array}{c}
                         4h_2(d_1d_3d_4) \\
                         2h_2(d_1d_3d_4) \\
                       \end{array}
                     \right) & \\\hline
a_5 &(0,1,0,1,1,0) & p_2p_3p_4\;\;\;\;\;\;p_2p_4\;\;\;\;\;\;p_3p_4 & \left(
                                                             \begin{array}{c}
                                                               -1 \\
                                                               +1 \\
                                                             \end{array}
                                                           \right)      & \la [\fp_2],[\fp_3]\ra & \\
& \la [\fp_2]\ra\;\; \la [\fp_3]\ra \;\;\la [\fp_2\fp_3]\ra & \la [\fp_3]\ra \;\;\; \la [\fp_2]\ra \;\;\;                                                                                                 \la [\fp_2],[\fp_3]\ra    & & D &32.033 \\
 & 2\qquad 2\qquad \left(
                     \begin{array}{c}
                       1\\
                       2 \\
                     \end{array}
                   \right)
  & 4\qquad 4\qquad \left(
                                       \begin{array}{c}
                                         h_2(d_1d_2) \\
                                         2h_2(d_1d_2) \\
                                       \end{array}
                                     \right)
&  & \left(
                                       \begin{array}{c}
                                         2h_2(d_1d_2) \\
                                         4h_2(d_1d_2) \\
                                       \end{array}
                                     \right) &  \\\hline
a_6 &(0,1,1,0,1,0) & p_1p_3p_4\;\;\;\;\;\;p_3p_4\;\;\;\;\;\;p_4 & \left(
                                                             \begin{array}{c}
                                                               -1 \\
                                                               +1 \\
                                                             \end{array}
                                                           \right)      & \la [\fp_1],[\fp_3]\ra & \\
& \la [\fp_1\fp_3]\ra\;\; \la [\fp_1]\ra \;\;\la [\fp_3]\ra & \la [\fp_1]\ra \;\;\; \la [\fp_1\fp_3]\ra \;\;\;                                                                                                 \la [\fp_1],[\fp_3]\ra    & & D & 32.033 \\
 & 2\qquad 2\qquad \left(
                     \begin{array}{c}
                       1\\
                       2 \\
                     \end{array}
                   \right)
  & 4\qquad 4\qquad \left(
                                       \begin{array}{c}
                                         h_2(d_1d_2) \\
                                         2h_2(d_1d_2) \\
                                       \end{array}
                                     \right)
&  & \left(
                                       \begin{array}{c}
                                         2h_2(d_1d_2) \\
                                         4h_2(d_1d_2) \\
                                       \end{array}
                                     \right) & \\\hline
a_7 &(0,1,0,0,1,0) & p_2p_4\;\;\;\;\;\;p_2p_4\;\;\;\;\;\; p_4 & \left(
                                                             \begin{array}{c}
                                                               -1 \\
                                                               +1 \\
                                                             \end{array}
                                                           \right)      & \la [\fp_1],[\fp_2]\ra  &\\
& \la [\fp_2]\ra\;\; \la [\fp_1]\ra \;\;\la [\fp_1\fp_2]\ra & \la [\fp_1]\ra \;\;\; \la [\fp_2]\ra \;\;\;                                                                                                 \la [\fp_1],[\fp_2]\ra    & & D & 64.144\\
 & 2\qquad 2\qquad \left(
                     \begin{array}{c}
                       1\\
                       2 \\
                     \end{array}
                   \right)
  & 4\qquad 4\qquad \left(
                                       \begin{array}{c}
                                         h_2(d_1d_2) \\
                                         2h_2(d_1d_2) \\
                                       \end{array}
                                     \right)
&  & \left(
                                       \begin{array}{c}
                                         2h_2(d_1d_2) \\
                                         4h_2(d_1d_2) \\
                                       \end{array}
                                     \right) &  \\\hline

a_8 &(1,0,1,0,1,0) & p_3p_4\;\;\; p_3p_4 \;\;\; \left(
                                                             \begin{array}{c}
                                                               p_1 \\
                                                               p_4 \;\text{or}\;p_1p_4 \\
                                                             \end{array}
                                                           \right)    & -1   & \la [\fp_1],[\fp_3]\ra &\\
& \la [\fp_3]\ra\;\; \la [\fp_1\fp_3]\ra \;\;\la [\fp_1]\ra & \la [\fp_1]\ra \;\;\; \left(
                                                                \begin{array}{c}
                                                                  \la [\fp_1]\ra \\
                                                                   \la [\fp_1],[\fp_3]\ra \\
                                                                \end{array}
                                                              \right)
 \;\;\; \la [\fp_3]\ra  & & \left(
                                                                                                             \begin{array}{c}
                                                                                                               S \\
                                                                                                               D \\
                                                                                                             \end{array}
                                                                                                           \right)
                                                                                   & 32.033             \\
 & 2\quad \left(
     \begin{array}{c}
       2 \\
       1 \\
     \end{array}
   \right)
   \quad 2 & 4\quad \left(
                       \begin{array}{c}
                         2h_2(d_1d_3d_4) \\
                         h_2(d_1d_3d_4) \\
                       \end{array}
                     \right)
   \quad 4  & & \left(
                       \begin{array}{c}
                         4h_2(d_1d_3d_4) \\
                         2h_2(d_1d_3d_4) \\
                       \end{array}
                     \right) & \\\hline
\rsp a_9 &(1,0,0,1,1,0) & p_1p_2p_4\;\;\;\;\;\;p_2p_4\;\;\;\;\;\;p_1p_4 & -1     & \la [\fp_1],[\fp_2]\ra & \\
& \la [\fp_1]\ra\;\; \la [\fp_2]\ra \;\;\la [\fp_1\fp_2]\ra & \la [\fp_1]\ra \;\;\; \la [\fp_2]\ra \;\;\; \la [\fp_1\fp_2]\ra  & & Q & 64.147 \\
 & 2\qquad 2\qquad 2 & 4\qquad 4\qquad 4 & & 8 & \\\hline
\rsp a_{10} &(1,1,1,0,0,0) & p_4\;\;\;\;\;\;p_4\;\;\;\;\;\; p_4 & -1     & \la [\fp_1],[\fp_2]\ra &  \\
& \la [\fp_1]\ra\;\; \la [\fp_2]\ra \;\;\la [\fp_1\fp_2]\ra & \la [\fp_1]\ra \;\;\; \la [\fp_2]\ra \;\;\; \la [\fp_1\fp_2]\ra  & & Q & 64.147 \\
 & 2\qquad 2\qquad 2 & 4\qquad 4\qquad 4 & & 8 &\\\hline
\rsp a_{11} &(1,1,0,1,1,0) & p_2p_3\;\;\;\;\;\;p_2p_4\;\;\;\;\;\; p_3p_4 & -1     & \la [\fp_1],[\fp_2]\ra & \\
& \la [\fp_1\fp_2]\ra\;\; \la [\fp_2]\ra \;\;\la [\fp_1]\ra & \la [\fp_1],[\fp_2]\ra \;\;\; \la [\fp_1],[\fp_2]\ra \;\;\; \la [\fp_1],[\fp_2]\ra  & & (2,2) & 32.39  \\
 & 1\qquad 1\qquad 1 & 2\qquad 2\qquad 2 & & 4 & \\\hline
\rsp a_{12} &(1,1,1,0,1,0) & p_2p_3\;\;\;\;\;\;p_3p_4\;\;\;\;\;\; p_4 & -1     & \la [\fp_1],[\fp_2]\ra & \\
& \la [\fp_1]\ra\;\; \la [\fp_1\fp_2]\ra \;\;\la [\fp_2]\ra & \la [\fp_1],[\fp_2]\ra \;\;\; \la [\fp_1],[\fp_2]\ra \;\;\; \la [\fp_1],[\fp_2]\ra  & & (2,2) & 32.39  \\
 & 1\qquad 1\qquad 1 & 2\qquad 2\qquad 2 & & 4 &\\\hline
\rsp a_{13} &(1,1,0,0,1,0) & p_2p_3\;\;\;\;\;\;p_2p_4\;\;\;\;\;\; p_4 & -1     & \la [\fp_1],[\fp_2]\ra & \\
& \la [\fp_1]\ra\;\; \la [\fp_2]\ra \;\;\la [\fp_1\fp_2]\ra & \la [\fp_1],[\fp_2]\ra \;\;\; \la [\fp_1],[\fp_2]\ra \;\;\; \la [\fp_1],[\fp_2]\ra  & & (2,2) & 32.34 \\
 & 1\qquad 1\qquad 1 & 2\qquad 2\qquad 2 & & 4 &\\\hline
\end{array}$$  }

{\tiny  $$
                            \right]   & -1    & \la [\fp_1],[\fp_2]\ra & \\
& \la [\fp_1]\ra\;\; \la [\fp_1\fp_2]\ra \;\;\la [\fp_2]\ra & \la [\fp_1],[\fp_2]\ra \; \;\;\la [\fp_1],[\fp_2]\ra \;\; \la [\fp_1],[\fp_2]\ra
 & & (2,2)
                                                                                        & 32.034     \\
 & 1\qquad 1
 \qquad 1  & 2\qquad 2
 \qquad 2  & & 4
   & \\\hline

\end{array}$$}

\newpage
\section{{\bf Some Properties of Relevant ${2}$-Groups}}\label{S6}
We now take a brief detour into the theory of finite $2$-groups that we will make use of in later sections. If we consider the groups $G^+$ with quotient groups $G^+/G^+_3$ of order $32$ in the sixth column of Appendix~II, we can see that these are contained in the family of all the finite  $2$-groups $G$ for which $G/G'\simeq (2,2,2)$ and $G'/G_3\simeq (2,2)$. (For $G^+/G^+_2\simeq \Cl_2^+(k)\simeq (2,2,2)$ and since $(G^+_2:G^+_3)=4$ and elementary (cf.\,\cite{Hup}, Satz 2.13), we have $G^+_2/G^+_3\simeq (2,2)$.) See below for more details. The groups $G$ may be given by $G=\la a_1,a_2,a_3 \ra$ where $G'=\la c_{12},c_{13},G_3\ra$ with $c_{23}\in G_3$. Here, $c_{ij}=[a_i,a_j]$ (and for later use $c_{ij\ell}=[a_i,a_j,a_\ell]=[c_{ij},a_\ell]$, etc.). Then the groups $G/G_3$ are of the following types with the given presentations and some properties of $G'$ (which we'll prove below). (The presentation and the structure of the Schur multiplier $\mathcal{M}$  of $G/G_3$ may be found in \cite{SW} with one change for this section only: we have permuted the indices $1$ and $3$ of the generators while leaving the index $2$ unchanged, to be consistent with some of our earlier work. Hence, for example, we have $c_{23}\in G_3$ instead of $c_{12}$.)

\begin{table}[ht]\label{Tab1}
$$ \begin{array}{c|ccc|c||c||c}
\text{Group Number} & a_1^2 & a_2^2 & a_3^2 & c_{23} & G' & \mathcal{M}(G/G_3) \\ \hline\hline
32.033  & 1 & 1 & 1 & 1 & d=2, 3 & (2,2,2,2)\\ \hline
32.034  & 1 & c_{12} & c_{13} & 1 & (2^*,2^*)& (2,2,4) \\ \hline
32.035  & c_{13} & c_{12} & c_{13} & 1 &(2,2^*) & (2,2)\\ \hline
32.036  & 1 & 1 & c_{13} & 1 & (2^*,2^*) & (2,2,2) \\ \hline
32.037  & c_{13} & 1 & c_{13} & 1& (2,2^*) & (2,2)\\ \hline
32.038  & 1 & c_{13} & 1 & 1 & (2,2^*) & (2,2)\\ \hline
32.039  & 1 & c_{12}c_{13} & c_{13} & 1& (2,2^*) & (2,4)\\ \hline
32.040  & c_{12} & c_{12}c_{13} & c_{13} & 1 & (2,2) & (2)\\ \hline
32.041  & 1 & c_{12}c_{13} & c_{12} & 1 &(2,2)& (2)\\ \hline
\end{array}$$
\caption{}
\end{table}

In this table, $d$ denotes rank, $(2,4)$ denotes the direct sum of cyclic groups of orders $2$ and $4$, etc., and $2^*$ denotes any power of $2$ greater than $1$. As an aside, observe the rather close relation between the possible $4$-rank of $G'$ and the (2-)rank of the corresponding Schur multiplier. The Group Number $32.0x$, for $x=33,\dots,41$, of $G/G_3$ denotes the $x$th group of order $32$ as listed in \cite{HS} (and \cite{SW}).
\bigskip

We'll start by recalling some ``commutator calculus'' and other relevant facts; see the exposition given, for example, in \cite{BS0}. Assume for the moment that $\fG$ is an  arbitrary group written multiplicatively. If $x,y\in \fG$, then define  $[x,y]=x^{-1}y^{-1}xy$, the commutator of $x$ with $y$. More generally, since $[*,*]$ is not associative, we define (inductively on $\ell$)
$[x_1,\dots,x_\ell]=[[x_1,\dots, x_{\ell-1}],x_\ell],$ for all $x_j\in \fG$ and $\ell >2$. If $A$ and $B$ are nonempty subsets of $\fG$, then $[A,B]=\la \{[a,b]:a\in A,b\in B\}\ra.$ In particular,  the commutator subgroup $\fG'$ is defined as $[\fG,\fG]$; also $\fG''=(\fG')'.$ Recall, too, that a group $\fG$ is metabelian if $\fG''=1$, i.e., $\fG'$ is abelian. Moreover, we denote
$$x^y=y^{-1}xy$$
for any $x,y\in\fG$. Notice then that
\begin{equation}\label{E0}
x^y=x[x,y]\,.
\end{equation}

Also notice that
\begin{equation}\label{E1}
[x,y]^z=[x^z,y^z]\,
\end{equation}
where $z\in\fG$.

Of particular use in our analysis is the lower central series $\{\fG_\ell\}$ of $\fG$, which is defined inductively as:
$\fG_1=\fG, \quad \text{and}\quad \fG_{\ell+1}=[\fG,\fG_\ell], \quad \text{for all}\;\; \ell\geq 1.$
In particular, notice that $\fG_2=\fG'.$  The lower central series is especially useful when the group $\fG$ is nilpotent, i.e., when
the lower central series terminates in finitely many steps at the identity subgroup $ 1$. As is well known, all finite $p$-groups are nilpotent, for  any prime $p$.

 Recall that (cf.\,\cite{Bla}) for any $x,y,z\in \fG$,
\begin{align}\label{E2}
    &[xy,z]=[x,z][x,z,y][y,z]\,,\\
    &[x,yz]=[x,z][x,y][x,y,z]\,,\notag\\
    &[x,y^{-1}]^y=[x^{-1},y]^x=[y,x]=[x,y]^{-1}\,,\notag\\
    &[x,y^{-1},z]^y[y,z^{-1},x]^z[z,x^{-1},y]^x=1\;\;\;\;(\mbox{{\em Witt identity}}), \notag\\
    & [x,y,z][y,z,x][z,x,y]\equiv 1 \bmod G''\;\;\;\;(\mbox{{\em Witt congruence}}). \notag
\end{align}

We now consider more specific groups. Let $G$ be a finite $2$-group such that $G/G'\simeq (2,2,2)$ and $G_2/G_3\simeq (2,2)$.   Hence,  $G/G_3$ is one of the groups of order $32$ given above, i.e., of type $32.0x$ for $x=33,\dots,41$.  Our first goal is to show that, with the exception of $G/G_3$ of type $32.033$,  all these groups $G$ have $G'$ of rank $2$. Moreover, in the exceptional case we will show that  the rank of $G'$ is at most $3$.

We start with a lemma:

\begin{lem}\label{L1}Let $G$ be a finite $2$-group such that $G/G_2\simeq (2,2,2)$ and $G_2/G_3\simeq (2,2)$, and  hence as given  in one of the rows of Table~1 above. Then

(a)\; $c_{ji\ell}\equiv c_{ij\ell}^{-1} \bmod G''$, for all $i,j,\ell\in \{1,2,3\}$;

(b)\; $c_{123}\equiv c_{132} \bmod G_4$.

\end{lem}
\begin{proof} In proving (a), notice that by (\ref{E2})
$$c_{ji\ell}=[c_{ij}^{-1},a_\ell]\equiv [c_{ij}^{-1},a_\ell]^{c_{ij}}=[a_\ell,c_{ij}]=c_{ij\ell}^{-1}\bmod G''.$$

For (b), notice by the Witt congruence that we have $1\equiv c_{123}c_{231}c_{312}\bmod G''$; thus also
$1\equiv c_{123}c_{231}c_{312}\bmod G_4$, since $G''\subseteq G_4$ for any $p$-group, see, e.g., Satz~2.12 of \cite{Hup}. But recall by Table~1 that  $c_{23}\in G_3$, and hence  $c_{231}=[c_{23},a_1]\in G_4$. Therefore, $c_{123}\equiv c_{312}^{-1}\equiv c_{132}\bmod G_4$, as we  wanted.

\end{proof}

Next, we prove a useful proposition.
\begin{prop}\label{P0} Let $G$ be a finite $2$-group such that $G/G_2\simeq (2,2,2)$ and $G_2/G_3\simeq (2,2)$, and hence   $G$ is given as in one of the rows of Table~1 above. Then for some positive integers $m,n$, we have

(a) \; $G'\simeq (2,2)$, if $G_3=G_4$;

(b) \;  $G'/G''\simeq (2,2^m)$, if $G_3=\la c_{12}^2,G_4\ra$;

(c) \;  $G'/G''\simeq (2^m,2^n)$, if $G_3=\la c_{12}^2,c_{13}^2, G_4\ra$.

\end{prop}
\begin{proof} With respect to (a), suppose $G_3=G_4$. Then since $G$ is nilpotent, $G_3=1$ and therefore $G'\simeq G_2/G_3\simeq (2,2).$

Now consider (b). Suppose $G_3=\la c_{12}^2, G_4\ra$. Since $G=\la a_1,a_2,a_3\ra$ and $G_2=\la c_{12},c_{13},c_{23},G_3\ra=\la c_{12},c_{13},G_3\ra$, we see that
$$G_3=\la c_{121},c_{122},c_{123},c_{131},c_{132},c_{133},G_4\ra.$$
Hence by part (b) of Lemma~\ref{L1}, we have
$$G_3=\la c_{121},c_{122},c_{123},c_{131},c_{133},G_4\ra.$$
By our hypothesis, for any $i,j,\ell\in \{1,2,3\}$,\; $c_{ij\ell}\equiv c_{12}^2$ or $c_{ij\ell}\equiv 1\bmod G_4$.
Now since $\exp(G/G')=2$, we know (see \cite{Hup}, Satz~2.13) that $\exp(G_j/G_{j+1})\leq 2$ for all positive integers $j$. Since $G_3=\la c_{12}^2,G_4\ra$,  $G_4=\la [c_{12}^2,a_1],[c_{12}^2,a_2],[c_{12}^2,a_3],G_5\ra.$  But by (\ref{E2}), we have $[c_{12}^2,a_j]=c_{12j}[c_{12j},c_{12}]c_{12j}\equiv c_{12j}^2\bmod G_5,$ since  $[c_{12j},c_{12}]\in [G_3,G_2]\subseteq G_5,$ cf., eg., \cite{Hup}, Hauptsatz~2.11. But then since $\exp(G_4/G_5)\leq 2$, we see from above that $c_{ij\ell}^2\equiv c_{12}^4$ or $c_{ij\ell}^2\equiv 1\bmod G_5$, and therefore $G_4=\la c_{12}^4, G_5\ra$. From this, we see that $G_3=\la c_{12}^2,G_5\ra$. Continuing (by induction) we have $G_3=\la c_{12}^2,G_j\ra$ for all $j\geq 3$. But then, since $G$ is nilpotent, $G_j=1$ for sufficiently large $j$. Therefore, $G_3=\la c_{12}^2\ra$ and so $G_2=\la c_{12},c_{13}\ra$. But this implies that $G'/G''\simeq (2,2^m)$ for some positive integer $m$. For $c_{13}^2=c_{12}^{2u}$ implies that $(c_{13}c_{12}^{-u})^2\equiv 1\bmod G''$, and thus $G_2=\la c_{12},c_{13}c_{12}^{-u}\ra.$

Finally, for (c), suppose $G_3=\la c_{12}^2,c_{13}^2,G_4\ra$. Then each $c_{ij\ell}\equiv c_{12}^{2a}c_{13}^{2b}\bmod G_4$ for some integers $a$ and $b$. But then $\la c_{12}^4,c_{13}^4,G_5\ra\subseteq G_4=\la c_{12i}^2,c_{13j}^2, G_5: i,j=1,2,3\ra\subseteq  \la c_{12}^4,c_{13}^4,G_5\ra$. Thus $G_4=\la c_{12}^4,c_{13}^4,G_5\ra$ and so $G_3=\la c_{12}^2,c_{13}^2,G_5\ra.$ Continuing, analogously to the previous case we see that  $G_3=\la c_{12}^2,c_{13}^2\ra$, whence $G_2=\la c_{12},c_{13}\ra$ and $G'/G''\simeq (2^m,2^n)$ for some positive integers $m$ and $n$. Notice  too (by induction) that $G_j=\la c_{12}^{2^{j-2}},c_{13}^{2^{j-2}}\ra=G_2^{2^{j-2}}$, for $j\geq 2$.
\end{proof}

From this proposition and some earlier results, we get the following corollary:

\begin{cor}\label{C1} Let $G$ be as in Proposition~\ref{P0} and satisfying any of the hypotheses in (a), (b), or (c). Then $G$ has derived length $2$, i.e., $G''=1$.
\end{cor}
\begin{proof} For (a) the result is trivial, since $G'\simeq (2,2)$ hence abelian. For (b), since  the $4$-rank of $G'/G''$ is at most $1$, we see that  $G''=1$ by \cite{Bla2}, Theorem~1. Finally for (c), the Main Theorem of \cite{BS2} shows that $G''=1$.
\end{proof}

We now want to prove the structure of $G'$ as given in Table~1 above. In light of the corollary, for all  but the first row in the table it suffices to prove the results for $G'/G''$.  We will consider all rows except the first one in that table. To this end, we need only compute $G_3/G_4$ by Proposition~\ref{P0}. For this, it suffices to determine  $[a_i^2,a_j]\bmod G_4$. First, notice that by (\ref{E2})
$$[a_i^2,a_j]\equiv c_{ij}^2c_{iji}\bmod G_4\,.$$
On the other hand, by the presentation modulo $G_3$ for each group, we will obtain a second condition on these commutator elements. For example, if $G/G_3$ is of type $32.041$, then $[a_2^2,a_1]\equiv [c_{12}c_{13},a_1]\equiv c_{121}c_{131}\bmod G_4.$

The following table gives these relations for each of the groups $32.0x$ for $x$ given in decreasing order: $41,\dots,34$.
Some of these group types and conclusions may be found in \cite{BLS2}. We have included them all here for the sake of completeness  and since it requires only  a small amount of space. In the table, the congruences should all be interpreted modulo $G_4$. Also recall that $c_{123}\equiv c_{132}\bmod G_4$ by Lemma~\ref{L1}.

\newpage

\begin{tiny}\begin{table}[ht]\label{Tab2}
$$ \begin{array}{c|ccc|c}
\text{Group Type} & [a_1^2,a_1] & [a_1^2,a_2] & [a_1^2,a_3]  & \mbox{Congruences} \\
  & [a_2^2,a_1] & [a_2^2,a_2] & [a_2^2,a_3]  &   \bmod G_4     \\
  & [a_3^2,a_1] & [a_3^2,a_2] & [a_3^2,a_3]  &  \mbox{Conclusions}     \\ \hline\hline
32.041 & 1\equiv 1 &1\equiv c_{12}^2c_{121}  & 1\equiv c_{13}^2c_{131} & 1\equiv c_{123}\equiv c_{122}\equiv c_{131}\equiv c_{13}^2    \\
  & c_{121}c_{131}\equiv c_{12}^2c_{122} & c_{122}c_{132}\equiv 1 &c_{123}c_{133}\equiv 1 & \equiv c_{133}\equiv c_{121} \\
  &  c_{121}\equiv c_{13}^2c_{133} & c_{122}\equiv 1 & c_{123}\equiv 1 &  G_3=G_4 \Rightarrow G_2\simeq (2,2) \\ \hline
32.040 & c_{121}\equiv 1 &c_{122}\equiv c_{12}^2c_{121}  & c_{123}\equiv c_{13}^2c_{131} & 1\equiv c_{121}\equiv c_{133}\equiv c_{131}\equiv c_{13}^2    \\
  & c_{121}c_{131}\equiv c_{12}^2c_{122} & c_{122}c_{132}\equiv 1 & c_{123}c_{133}\equiv 1 & \equiv c_{123}\equiv c_{122} \\
  &  c_{131}\equiv c_{13}^2c_{133} & c_{132}\equiv 1 & c_{133}\equiv 1 &  G_3=G_4 \Rightarrow G_2\simeq (2,2) \\ \hline
32.039 & 1\equiv 1 &1\equiv c_{12}^2c_{121}  & 1\equiv c_{13}^2c_{131} & 1\equiv c_{133}\equiv c_{122}\equiv c_{123}\equiv c_{131},    \\
  & c_{121}c_{131}\equiv c_{12}^2c_{122} & c_{122}c_{132}\equiv 1 &c_{123}c_{133}\equiv 1 &  c_{121}\equiv c_{12}^2 \\
  &  c_{131}\equiv c_{13}^2c_{133} & 1\equiv 1 & c_{133}\equiv 1 &  G_3=\la c_{12}^2,G_4\ra \Rightarrow G'/G''\simeq (2,2^*) \\ \hline
32.038 & 1\equiv 1 &1\equiv c_{12}^2c_{121}  & 1\equiv c_{13}^2c_{131} & 1\equiv c_{123}\equiv c_{133}\equiv c_{13}^2\equiv c_{131},    \\
  & c_{131}\equiv c_{12}^2c_{122} & c_{132}\equiv 1 & c_{133}\equiv 1 &  c_{121}\equiv c_{122}\equiv c_{12}^2 \\
  &  1\equiv c_{13}^2c_{133} & 1\equiv 1 & 1\equiv 1 &  G_3=\la c_{12}^2,G_4\ra \Rightarrow G'/G''\simeq (2,2^*) \\ \hline
32.037 & c_{131}\equiv 1 & c_{132}\equiv c_{12}^2c_{121}  & c_{133}\equiv c_{13}^2c_{131} & 1\equiv c_{133}\equiv c_{123}\equiv  c_{131},    \\
  & 1\equiv c_{12}^2c_{122} & 1\equiv 1 & 1 \equiv 1 &  1\equiv c_{122}\equiv c_{121}\equiv c_{12}^2 \\
  & c_{131}\equiv c_{13}^2c_{133} & c_{132}\equiv 1 & c_{133}\equiv 1 &  G_3=\la c_{12}^2,G_4\ra \Rightarrow G'/G''\simeq (2,2^*) \\ \hline
32.036 & 1\equiv 1 & 1\equiv c_{12}^2c_{121}  & 1\equiv c_{13}^2c_{131} & 1\equiv c_{133}\equiv c_{123},   \\
  & 1\equiv c_{12}^2c_{122} & 1\equiv 1 & 1\equiv 1 &  c_{121}\equiv c_{12}^2\equiv c_{122},\; c_{131}\equiv c_{13}^2 \\
  &  c_{131}\equiv c_{13}^2c_{133} & c_{132}\equiv 1 & c_{133}\equiv 1 &  G_3=\la c_{12}^2,c_{13}^2,G_4\ra \Rightarrow G'/G''\simeq (2^*,2^*) \\ \hline
32.035 & c_{131}\equiv 1 & c_{132}\equiv c_{12}^2c_{121}  & c_{133}\equiv c_{13}^2c_{131} & 1\equiv c_{133}\equiv c_{123}\equiv  c_{131} \equiv c_{122},   \\
  & c_{121}\equiv c_{12}^2c_{122} & c_{122}\equiv 1 & c_{123}\equiv 1 &  1\equiv c_{121}\equiv c_{12}^2 \\
  & c_{131}\equiv c_{13}^2c_{133} & c_{132}\equiv 1 & c_{133}\equiv 1 &  G_3=\la c_{12}^2,G_4\ra \Rightarrow G'/G''\simeq (2,2^*) \\ \hline
32.034 & 1\equiv 1 & 1\equiv c_{12}^2c_{121}  & 1\equiv c_{13}^2c_{131} & 1\equiv c_{133}\equiv c_{123}\equiv c_{122},   \\
  & c_{121}\equiv c_{12}^2c_{122} & c_{122}\equiv 1 & c_{123}\equiv 1 &  c_{121}\equiv c_{12}^2,\; c_{131}\equiv c_{13}^2 \\
  &  c_{131}\equiv c_{13}^2c_{133} & c_{132}\equiv 1 & c_{133}\equiv 1 &  G_3=\la c_{12}^2,c_{13}^2,G_4\ra \Rightarrow G'/G''\simeq (2^*,2^*) \\ \hline
\end{array}$$
\caption{}
\end{table} \end{tiny}

%\newpage

Now we see that $G'$ has rank $2$ for all these groups and that $G''=1$. This is not the case for $G/G_3$ of type $32.033$, however, as we have seen in \cite{BLS2}. We do, nevertheless, have a bound on the rank of $G'$ in this case:

\begin{prop}\label{P01} Let $G$ be a finite $2$-group for which $G/G_3$ is of type $32.033$. Then $\rank(G')\leq 3$.
\end{prop}
\begin{proof} Without loss of generality, we assume that $$G_2^2=1\,, $$ by the Burnside Basis Theorem, see \cite{Hup}. Thus we also have $G''=1$, since $G''\subseteq G_2^2$ (Satz 3.14 of \cite{Hup}, for example). Using the presentation modulo $G_3$ as in Table~1, we see that $G=\la a_1,a_2,a_3\ra$ with $a_1^2\equiv a_2^2\equiv a_3^2\equiv c_{23}\equiv 1\bmod G_3.$ Thus $G_2=\la c_{12},c_{13}, G_3\ra$ (as we have already seen).

We will show that $G_3=\la c_{123}, G_4\ra $ which implies that $G_2=\la c_{12},c_{13},c_{123},G_4\ra$. We will then show that $G_4=1$, which will establish the proposition.

First, as noted before, $G_3=\la c_{121},c_{122},c_{131},c_{133},c_{123},G_4\ra.$ Now since $a_j^2\equiv 1~\bmod G_3$, we have, by (\ref{E2}) and our assumption that $G_2^2=1$, that $1\equiv [a_i,a_j^2]=c_{ij}^2c_{ijj}\equiv c_{ijj}\bmod G_4.$  Thus by Lemma~\ref{L1}, $G_3=\la c_{123},G_4\ra$. From this we have $$G_4=\la c_{1231},c_{1232},c_{1233},G_5\ra\,.$$

Next, since $G''=1$, the Witt congruence in (\ref{E2}) becomes an equality. In particular, we have
$$1=[c_{12},a_3,a_1][a_3,a_1,c_{12}][a_1,c_{12},a_3]=c_{1231}[c_{31},c_{12}][c_{121}^{-1},a_3]\equiv c_{1231}\bmod G_5\,,$$
since $[c_{31},c_{12}]=1$ and $c_{121}\equiv 1\bmod G_4$ by our presentation of $G/G_3$. Similarly, we have $c_{1232}\equiv 1\bmod G_5.$ With respect to $c_{1233}$, it will be more convenient to replace this element with another. Namely, notice that since $c_{123}\equiv c_{132}\bmod G_4$, by Lemma~\ref{L1}, $c_{1233}\equiv c_{1323}\bmod G_5$ and therefore we may replace $c_{1233}$ with $c_{1323}$.
But now proceeding as before, we find that $c_{1323}\equiv 1\bmod G_5$. Therefore $G_4=G_5$ and so $G_4=1$, as desired.
\end{proof}

\section{{\bf The Structure of ${G^+/G_3^+}$ for ${G^+=\Gal(k_+^2/k)}$}}\label{S7}

As before, we consider real quadratic fields $k$ for which $\Cl_2(k)\simeq (2,2)$ and where $d_k$ is not a sum of two squares. We are interested in determining properties of the narrow $2$-class field tower over $k$. Toward this end, we will start by studying $G^+=\Gal(k_+^2/k)$, or more accurately, $G^+/G_3^+$.
Since $\Cl_2(k)\simeq (2,2)$ and $d_k$ is not the sum of two squares, the norm of the fundamental unit of $k$ satisfies $N\eps_k=+1$ and $\Cl_2^+(k)\simeq (2,2,2)$.
Hence $d_k=d_1d_2d_3d_4$ is in one of the Types I, II, III, IV.  This all implies that the Hilbert $2$-class field $k^1$ of $k$ coincides with the genus field, which equals $k_\gen=\Q(\sqrt{d_1d_4},\sqrt{d_2d_4},\sqrt{d_3d_4}\,)$ when all four prime discriminants are negative, and $\Q(\sqrt{d_1},\sqrt{d_2},\sqrt{d_3d_4}\,)$, if not. Moreover, the ordinary $2$-class field tower has length at most $2$. The Galois group of the second $2$-class field $k^2$ over $k$ is also well known and is equal to $\Gal(k^1/k)=\Cl_2(k)$ for  $d_k$ of Types III, IV and as given in Appendix II above for Types I, II.

We now consider $k_+^1=k^+_\gen=\Q(\sqrt{d_1}, \sqrt{d_2},\sqrt{d_3},\sqrt{d_4}\,),$ the strict genus field of $k$ and the Galois group $\Gal(k_+^1/k)\simeq (2,2,2)$, and will study the structure of  the $2$-class field tower over $k_+^1$. Notice that since $k_+^1$ is totally imaginary, its narrow and wide $2$-class field towers coincide.

In Appendix III below, we  give  a presentation of $G^+/G_3^+$ in terms of quadratic residue symbols according to Koch's description given in \cite{Koc}. First some notation: For  $i,j\in\{1,2,3,4\}$  let  $\nu_{ij}=0$ or $1$ according as
$$(-1)^{\nu_{ij}}=\bigg(\frac{d_i}{p_j}\bigg)$$
where $p_j$ is the prime dividing $d_j$ and $(\cdot/\cdot)$ is the Kronecker symbol when $i\not=j$ and $(d_i/p_i)=(d_i'/p_i)$ such that $d_i'=d_k/d_i$.

We write $\mu_i$ for $\nu_{i4}$ and $\delta_i$ for $\nu_{ii}$. (This is not to be confused with the $\delta$'s in the tables in Appendix II.)

For the $k$ satisfying the conditions above, we have $G^+/G_3^+=\la s_1,s_2,s_3\ra$ such that
$$\begin{array}{lcl}
  s_1^{2\delta_1}s_2^{2\nu_{21}}s_3^{2\nu_{31}} & = & t_{12}^{\nu_{21}}t_{13}^{\nu_{31}} \\
  s_1^{2\nu_{12}}s_2^{2\delta_2}s_3^{2\nu_{32}} & = & t_{12}^{\nu_{12}}t_{23}^{\nu_{32}} \\
  s_1^{2\nu_{13}}s_2^{2\nu_{23}}s_3^{2\delta_3} & = & t_{13}^{\nu_{13}}t_{23}^{\nu_{23}} \\
  s_1^{2\mu_{1}}s_2^{2\mu_{2}}s_3^{2\mu_{3}} & = & 1 \\
  G^+_3& = & 1\,,
\end{array}$$
where $t_{ij}=[s_i,s_j]$.  We present $G^+/G^+_3$  for the four cases described above.
For  each type, we   give first the $9$-tuple representing the vector \\ $(\nu_{12},\nu_{13},\nu_{23};\mu_1,\mu_2,\mu_3;\delta_1,\delta_2,\delta_3);$
secondly, the Koch presentation involving $s_1,s_2,s_3$; thirdly, a transformation (if necessary) of the $s_i$ into a new presentation involving $a_1,a_2,a_3$; and then finally the presentation involving the $a_j$ and its isomorphism type as given in \cite{HS} or \cite{SW}. We let $c_{ij}=[a_i,a_j]$ below.

In this appendix, we note that in the data for Types I and  II, $\nu_{ij}=\nu_{ji}$ whereas for Types III and IV, $\nu_{ij}=1-\nu_{ji},$  \;($i\not=j$). Moreover, with respect to the data of Type II, we have only included the cases with $\mu_3=0$, since the data for $\mu_3=1$ are identical, because $s_3^2=1.$

\section{{\bf Appendix III}}

$$\begin{array}{l|l|l}
{\rm Type~I} & (\nu_{12},\nu_{13},\nu_{23};\mu_1,\mu_2,\mu_3;\delta_1,\delta_2,\delta_3) & G^+/G_3^+\\\hline
(a_1) & (0,0,0;1,1,0;1,1,1) & \\
      & s_1^2=s_2^2=s_3^2=1 &  64.144 \\\hline
(a_2) & (0,1,1;0,0,0;1,1,1) & \\
      & s_1^2=t_{13}, s_2^2=t_{23}, s_3^2=1 &  \\
      & a_1=s_1s_3, a_2=s_2s_3, a_3=s_3 & \\
      & a_1^2=a_2^2=a_3^2=1 & 64.144 \\\hline
(a_3) & (1,0,1;0,0,0;1,0,0) & \\
      & s_1^2=t_{12}t_{23},s_2^2=t_{23}, s_3^2=1  &  \\
      & a_1=s_1s_2, a_2=s_2s_3, a_3=s_3 &\\
      & a_1^2=a_2^2=a_3^2=1 & 64.144\\\hline
(a_4) & (1,0,0;0,1,0;1,0,1) & \\
      & s_1^2=t_{12},s_2^2=s_3^2=1  &  \\
      & a_1=s_2, a_2=s_1s_2, a_3=s_3 &\\
      & a_1^2=a_2^2=a_3^2=1 & 64.144\\\hline
(a_5) & (0,1,0;1,1,0;0,1,0) & \\
      & s_1^2=s_2^2=s_3^2=1;t_{13}=1  &  \\
      & a_1=s_1, a_2=s_3, a_3=s_2 &\\
      & a_1^2=a_2^2=a_3^2=1; c_{12}=1 & 32.033\\\hline
(a_6) & (0,1,1;0,1,0;1,0,1) & \\
      & s_1^2=t_{13},s_2^2=s_3^2=1;t_{23}=1  &  \\
      & a_1=s_3, a_2=s_2, a_3=s_1s_3 &\\
      & a_1^2=a_2^2=a_3^2=1; c_{12}=1 & 32.033\\\hline
(a_7) & (0,1,0;0,1,0;1,1,0) & \\
      & s_1^2=t_{13},s_2^2=s_3^2=1  &  \\
      & a_1=s_2, a_2=s_3, a_3=s_1s_3 &\\
      & a_1^2=a_2^2=a_3^2=1 & 64.144\\\hline
(a_8) & (1,0,1;0,1,0;1,1,0) & \\
      & s_1^2=t_{12}, s_2^2=s_3^2=1; t_{23}=1  &  \\
      & a_1=s_2, a_2=s_3, a_3=s_1s_2 &\\
      & a_1^2=a_2^2=a_3^2=1; c_{12}=1 & 32.033\\\hline
(a_9) & (1,0,0;1,1,0;0,0,1) & \\
      & s_1^2=s_2^2=t_{12}, s_3^2=1 &  64.147 \\\hline
(a_{10}) & (1,1,1;0,0,0;0,0,1) & \\
      & s_1^2=t_{12}t_{23},s_2^2=t_{12}t_{13},s_3^2=1  &  \\
      & a_1=s_1s_3, a_2=s_2s_3, a_3=s_3 &\\
      & a_1^2=a_2^2=c_{12}, a_3^2=1 & 64.147\\\hline
(a_{11}) & (1,1,0;1,1,0;1,0,0) & \\
      & s_1^2=s_2^2=t_{12}, s_3^2=1; t_{12}t_{13}=1  &  \\
      & a_1=s_1, a_2=s_1s_2s_3, a_3=s_3 &\\
      & a_1^2=c_{13}, a_2^2=c_{13}c_{23}, a_3^2=1; c_{12}=1 & 32.039\\\hline
(a_{12}) & (1,1,1;0,1,0;0,1,1) & \\
      & s_1^2=t_{13}t_{23},s_2^2=s_3^2=1; t_{12}t_{13}=1  &  \\
      & a_1=s_2s_3, a_2=s_1, a_3=s_2 &\\
      & a_1^2=c_{13}, a_2^2=c_{13}c_{23}, a_3^2=1; c_{12}=1 & 32.039\\\hline
(a_{13}) & (1,1,0;0,1,0;0,0,0) & \\
      & s_1^2=t_{13},s_2^2=s_3^2=1; t_{12}t_{13}=1  &  \\
      & a_1=s_2s_3, a_2=s_1, a_3=s_3 &\\
      & a_1^2=c_{13}, a_2^2=c_{23}, a_3^2=1; c_{12}=1 & 32.034\\\hline
\end{array}$$

$$\begin{array}{l|l|l}
{\rm Type~II} & (\nu_{12},\nu_{13},\nu_{23};\mu_1,\mu_2,\mu_3;\delta_1,\delta_2,\delta_3) & G^+/G_3^+\\\hline
(b_1) & (0,1,1;1,1,0;1,1,1) & \\
      & s_1^2=s_2^2=t_{13}, s_3^2=1; t_{13}t_{23}=1 & \\
      & a_1=s_1s_2, a_2=s_3,a_3=s_1s_3 & \\
      & a_1^2=c_{13}, a_2^2=a_3^2=1;c_{12}=1 & 32.036 \\\hline
(b_2) & (0,1,1;0,0,0;1,1,1) & \\
      & s_1^2=t_{13}, s_2^2=t_{23},s_3^2=1 & \\
      & \mbox{see}\;\;(a_2) & 64.144 \\\hline
(b_3) & (1,0,1;0,0,0;1,0,0) & \\
      & s_1^2=t_{12}t_{23}, s_2^2=t_{23}, s_3^2=1 & \\
      & \mbox{see}\;\;(a_3) & 64.144 \\\hline
(b_4) & (1,0,1;0,1,0;1,0,0) & \\
      & s_1^2=t_{12}, s_2^2=s_3^2=1; t_{23}=1 & \\
      & \mbox{see}\;\;(a_8) & 32.033 \\\hline
(b_5) & (0,1,0;1,1,0;1,0,0) & \\
      & s_1^2=s_2^2=t_{13}, s_3^2=1 & \\
      & a_1=s_3, a_2=s_1s_3,a_3=s_1s_2 & \\
      & a_1^2=a_2^2=1, a_3^2=c_{13}c_{23} & 64.146 \\\hline
(b_6) & (0,1,1;0,1,0;1,1,1) & \\
      & s_1^2=t_{13}, s_2^2=s_3^2=1; t_{23}=1 & \\
      & \mbox{see}\;\;(a_6) & 32.033 \\\hline
(b_7) & (0,1,0;0,1,0;1,0,0) & \\
      & s_1^2=t_{13}, s_2^2=s_3^2=1 & \\
      &\mbox{see}\;\;(a_7) & 64.144 \\\hline
(b_8) & (1,0,0;0,1,0;1,1,1) & \\
      & s_1^2=t_{12}, s_2^2=s_3^2=1 & \\
      & \mbox{see}\;\;(a_4) & 64.144 \\\hline
(b_9) & (1,1,1;1,1,0;0,0,1) & \\
      & s_1^2=s_2^2=t_{12}t_{13}, s_3^2=1; t_{13}t_{23}=1 & \\
      & a_1=s_1s_2, a_2=s_3,a_3=s_2s_3 & \\
      & a_1^2=c_{13}, a_2^2=1, a_3^2=c_{13};c_{12}=1 & 32.037 \\\hline
(b_{10}) & (1,1,1;0,0,0;0,0,1) & \\
      & s_1^2=t_{12}t_{23}, s_2^2=t_{12}t_{13}, s_3^2=1 & \\
      & \mbox{see}\;\;(a_{10}) & 64.147 \\\hline
(b_{11}) & (1,1,0;1,1,0;0,1,0) & \\
      & s_1^2=s_2^2=t_{13}, s_3^2=1; t_{12}=1 & \\
      & a_1=s_1, a_2=s_1s_2,a_3=s_3 & \\
      & a_1^2=c_{13}, a_2^2=a_3^2=1;c_{12}=1 & 32.036 \\\hline
(b_{12}) & (1,1,1;0,1,0;0,0,1) & \\
      & s_1^2=t_{13}t_{23},s_2^2=s_3^2=1; t_{12}t_{13}=1 & \\
      & \mbox{see}\;\;(a_{12})& 32.039 \\\hline
(b_{13}) & (1,1,0;0,1,0;0,1,0) & \\
      & s_1^2=t_{13},s_2^2=s_3^2=1; t_{12}t_{13}=1 & \\
      & \mbox{see}\;\;(a_{13})  & 32.034 \\\hline
\end{array}$$

\newpage

\bigskip

$$\begin{array}{l|l|l}
{\rm Type~III} & (\nu_{12},\nu_{13},\nu_{23};\mu_1,\mu_2,\mu_3;\delta_1,\delta_2,\delta_3) & G^+/G_3^+\\\hline
(c_1) & (1,0,1;1,0,0;1,0,0) & \\
      & s_1^2=1, s_2^2=t_{23}, s_3^2=t_{13}; t_{12}=1 & \\
      & a_1=s_2, a_2=s_1,a_3=s_1s_3 & \\
      & a_1^2=c_{13}, a_2^2=a_3^2=1;c_{12}=1 & 32.036 \\\hline
(c_2) & (1,0,0;1,0,0;1,1,1) & \\
      & s_1^2=1, s_2^2=t_{12}t_{23}, s_3^2=1; t_{13}=1 & \\
      & a_1=s_1, a_2=s_1s_3,a_3=s_1s_2s_3 & \\
      & a_1^2=a_2^2=a_3^2=1;c_{12}=1 & 32.033 \\\hline
(c_3) & (1,1,0;1,1,0;0,0,0) & \\
      & s_1^2=s_2^2=t_{13}, s_3^2=t_{12}t_{13}t_{23} & \\
      & a_1=s_1, a_2=s_1s_3,a_3=s_1s_2 & \\
      & a_1^2= c_{12},a_2^2=c_{23},a_3^2=c_{13} & 64.150 \\\hline
\end{array} $$

\bigskip
\bigskip
\bigskip

$$\begin{array}{l|l|l}
{\rm Type~IV} & (\nu_{12},\nu_{13},\nu_{23};\mu_1,\mu_2,\mu_3;\delta_1,\delta_2,\delta_3) & G^+/G_3^+\\\hline
(d_1) & (1,0,1;0,0,0;0,0,0) & \\
      & s_1^2=t_{12}, s_2^2=t_{23}, s_3^2=t_{13} & 64.150\\ \hline
(d_2) & (1,0,1;0,1,0;0,0,0) & \\
      & s_1^2=t_{12}, s_2^2=1, s_3^2=t_{13}; t_{23}=1 &  \\
      & a_1=s_3, a_2=s_2,a_3=s_1s_2 & \\
      & a_1^2=c_{13}, a_2^2=a_3^2=1;c_{12}=1 & 32.036 \\\hline
(d_3) & (1,0,1;1,1,0;0,0,0) & \\
      & s_1^2=s_2^2=t_{12},s_3^2=t_{13}; t_{12}t_{23}=1 &  \\
      & a_1=s_2, a_2=s_1s_2s_3,a_3=s_1 & \\
      & a_1^2=c_{13}, a_2^2=1, a_3^2=c_{13};c_{12}=1 & 32.037 \\\hline
(d_4) & (1,0,1;1,1,1;0,0,0) & \\
      & s_1^2=t_{12}, s_2^2=t_{23}, s_3=t_{13}; t_{12}t_{13}t_{23}=1 &  \\
      & a_1=s_1s_3, a_2=s_1s_2,a_3=s_1s_2s_3 & \\
      & a_1^2=c_{23}, a_2^2=c_{13}c_{23}, a_3^2=1;c_{12}=1 & 32.041 \\\hline
(d_5) & (1,1,1;0,1,0;1,0,1) & \\
       & s_1^2=s_2^2=1, s_3^2=t_{13}t_{23}; t_{12}=1 &  \\
       & a_1=s_1, a_2=s_2,a_3=s_1s_2s_3 & \\
      & a_1^2=a_2^2=a_3^2=1;c_{12}=1 & 32.033 \\\hline
(d_6) & (1,1,1;0,0,1;1,0,1) & \\
      & s_1^2=1, s_2^2=t_{13}t_{23},s_3^2=1; t_{12}=1 &  \\
      & a_1=s_1s_2, a_2=s_1,a_3=s_1s_2s_3 & \\
      & a_1^2=c_{13}, a_2^2=1, a_3^2=1;c_{12}=1 & 32.036 \\\hline
(d_7) & (1,1,1;1,0,1;1,0,1) & \\
      & s_1^2=1, s_2^2=t_{13}t_{23},s_3^2=1; t_{12}=1 &  \\
       & \mbox{see ($d_6$)} & 32.036 \\\hline
(d_8) & (1,1,1;1,1,0;1,0,1) & \\
      & s_1^2=s_2^2=1, s_3^2=t_{13}t_{23}; t_{12}=1 &  \\
       & \mbox{see ($d_5$)} & 32.033 \\\hline
\end{array} $$

\newpage

\section{{\bf The Rank of $\Cl_2(k_+^1)$}}\label{S8}

Given the information on $G^+/G_3^+$ in Appendix III, we are able to deduce the following results about the rank of $\Cl_2(k_+^1)$, the $2$-class group of the narrow $2$-class field of $k$. As before, we let $G^+=\Gal(k_+^2/k)$.

\begin{thm}\label{T1} Let $k$ be a real quadratic number field for which $\Cl_2(k)\simeq (2,2)$ and for which the discriminant $d_k$ is not the sum of two squares. Then
$$\rank\Cl_2(k_+^1)\left\{
                     \begin{array}{ll}
                       \geq 3, &  \hbox{if} \;\;(G^+:G_3^+)=64,\\
                       =3, & \hbox{if} \;\;G^+/G_3^+\simeq 32.033,\\
                       =2, & \hbox{otherwise}.
                     \end{array}
                   \right.$$
\end{thm}
(Observe that by Appendices I and III, the conditions in this proposition are determined by relations among the quadratic residue symbols involving the prime discriminants dividing $d_k$. As an example of this, see the next proposition for the determination of those $k$ for which $G^+/G_3^+\simeq 32.033$.)

\begin{proof} Suppose $G^+/G^+_3$ has order $64$. Since $G^+/G^+_2\simeq (2,2,2)$, hence elementary,  then $G^+_2/G^+_3$ is elementary and of order $8$. Thus $G^+_2/G^+_3\simeq (2,2,2)$. Hence $\rank\Cl_2(k_+^1)=\rank G_2^+\geq 3.$
When $G^+/G^+_3$ has order $32$ and not of type $32.033$, then the result follows immediately from Table~1 in Section~6.

Now suppose $G^+/G_3^+$ is of type $32.033$.  Then as an application of Scholz's knot sequence, see, for example, Section~3 of \cite{BLS0}, we have
$$\rank\Cl_2(k_+^1)=\rank\Gal(k_+^2/k_+^1)=\rank(G_2^+)\geq \rank\Gal(K_\cen^+/K_\gen^+)$$ $$\geq \rank M(G^+/G_3^+)-\rank E_k^+/E_k^+\cap N_{K/k}K^\times,$$
where $K$ is the fixed field of $G_3^+$ in $k_+^2$; $K_\cen^+$ and $K_\gen^+$ are the narrow $2$-central and $2$-genus class fields of $K/k$; and $E_k^+$ is the subgroup of totally positive units in $E_k$. But since $M(G^+/G_3^+)\simeq (2,2,2,2)$ (hence of rank $4$ for $G^+/G_3^+$ of type $32.033$) and since $\rank E_k^+=1$, we see that $\rank \Cl_2(k_+^1)\geq 3$ and therefore $=3$, by Proposition~\ref{P01}.
\end{proof}

(The inequalities are not sharp. For example, if $G^+/G_3^+\simeq 64.144$, then the Schur multiplier is $(2,2,2,2,2)$ and thus the lower bound on the rank is at least $4$.
Moreover, with the exception of the case $64.150$, in all the  cases for which $G^+/G_3^+$ has order $64$, all our examples have  rank  at least  $4$. It turns out that  if $G$ is any finite $2$-group for which $G/G_3\simeq 64.150$, then $\rank G_2=3$,  as we'll see in Part II.)

\medskip

As promised above, we now give a characterization of those $k$ for which $G^+/G_3^+\simeq 32.033$.

\begin{prop}\label{P3} Let $k$ satisfy the conditions in Theorem~\ref{T1}. Then $G^+/G_3^+\simeq 32.033$ if and only if $d_k=d_1d_2d_3d_4$ for some prime discriminants and after possible reordering of the $d_i$, one of the following conditions holds:

(1)  $d_1,d_2>0, d_3,d_4<0$   such that  $d_4=-4$ if $d_k\equiv 4\bmod 8$, and
$$(d_1/p_4)=+1, (d_2/p_3)=(d_2/p_4)=(d_1/p_2p_3)=-1,$$

(2) $d_i<0$ for all $i$ and
 $$ (d_1/p_2)=(d_1/p_3)=(d_2/p_3)=(d_2/p_4)=(d_4/p_1)=(d_4/p_3)=-1,$$
 and if $d_4=-4$, then in addition $(d_3/p_4)=+1.$
\end{prop}

This is immediate since by Appendix III, $G^+/G_3^+\simeq 32.033$ if and only if $d_k$ satisfies one of the cases $(a_5), (a_6), (a_8),(b_4), (b_6), (c_2), (d_5), (d_8)$ as given in Appendix I or I$'$.

\bigskip
We now give examples, using Pari, of those $k$ in each of the cases $(a_i),\dots, (d_i)$. In the table below we display the type of group $G^+/G^+_3$, the relevant case, the discriminant $d_k$ of the field $k$, and $\Cl_2(k_+^1)$.

\begin{tiny}\begin{table}[ht]
$$\begin{array}{l|l|l|l}
\rsp \text{Group type} & {\rm Case}  & d_k= & \Cl_2(k_+^1)\simeq\\\hline
\rsp 32.033 & (a_5) & 89\cdot 5\cdot 19\cdot 7 & (2,4,8)\\
 \rsp & (a_6) & 5\cdot 89\cdot 3\cdot 19   & (2,4,8)\\
\rsp & (a_8) & 13\cdot5\cdot 23\cdot3 & (2,4,8)\\
\rsp & (b_4) & 17\cdot 29\cdot 19\cdot 4 & (4,4,4)\\
\rsp & (b_6) & 17\cdot 13\cdot 7\cdot 4  & (2,4,8)\\
\rsp & (c_2) & 7\cdot 3\cdot 43\cdot 31  & (2,8,8)\\
\rsp & (d_5) & 31\cdot 3\cdot 23\cdot 4  & (2,4,4)\\
\rsp & (d_8) & 3\cdot 11\cdot 263\cdot 4 & (2,8,16)\\\hline
\rsp 32.034 & (a_{13}) & 13\cdot 5\cdot 11\cdot 23 & (4,8) \\
\rsp       & (b_{13}) & 17\cdot 5\cdot 11\cdot 4 & (4,4)  \\\hline
\rsp 32.036 & (b_1) & 29\cdot 13\cdot 11\cdot 4   & (4,8)  \\
\rsp       &(b_{11}) & 13\cdot 5\cdot 11\cdot 4  &  (4,4) \\
\rsp       & (c_1)  &7\cdot 3\cdot 11\cdot 31  & (4,8) \\
\rsp      & (d_2)  &  7\cdot 3\cdot 23\cdot 4  &  (4,8) \\
\rsp      & (d_6)  & 23\cdot 7\cdot 19\cdot 4  & (4,8)   \\
\rsp       & (d_7) & 3\cdot 23\cdot 11 \cdot 4 & (4,4)  \\\hline
\rsp 32.037 & (b_9) &  13\cdot 5\cdot 7\cdot 4 & (2,8)  \\
\rsp       & (d_3) & 3\cdot 11\cdot 79\cdot 4  & (2,16)  \\\hline
\rsp 32.039 & (a_{11}) & 13\cdot 5\cdot 11\cdot 67 & (2,4) \\
\rsp      & (a_{12}) & 5\cdot 13\cdot 7\cdot 11  & (2,4)\\
\rsp       & (b_{12}) & 17\cdot 5\cdot 7\cdot 4 & (2,4)  \\\hline
\rsp 32.041 & (d_4) & 3\cdot 11\cdot 19\cdot 4 & (2,2) \\\hline
\rsp 64.144 & (a_1) & 17\cdot 13\cdot 43\cdot 31 & (2,2,2,4,16)  \\
 \rsp      & (a_2) & 13\cdot 17\cdot 7\cdot 43 & (2,4,4,16) \\
\rsp       & (a_3) & 5\cdot 17\cdot 31\cdot 59 &  (2,2,4,4,8)\\
\rsp       & (a_4) & 5\cdot 17\cdot 19 \cdot 11  & (2,4,4,8)  \\
 \rsp      & (a_7) & 17\cdot 13\cdot 3\cdot 19 & (2,4,4,8) \\
\rsp       & (b_2) & 17\cdot 89 \cdot 7\cdot 4 & (2,2,2,4,4) \\
\rsp       & (b_3) & 97\cdot 17 \cdot 3\cdot 4  & (2,4,8,8) \\
 \rsp      & (b_7) &  17\cdot 13\cdot 3\cdot 4 &  (2,4,8,8) \\
\rsp       & (b_8) & 17\cdot 5\cdot 19\cdot 4 & (2,2,2,4,8) \\\hline
\rsp 64.146 & (b_5) & 13\cdot 29\cdot 7\cdot 4  & (2,2,8,32) \\\hline
\rsp 64.147 & (a_9) & 13\cdot 5\cdot 131\cdot 7 & (2,2,2,2,8) \\
 \rsp      & (a_{10}) & 17\cdot 5\cdot 3\cdot 19  &  (2,2,2,2,8)\\
 \rsp      & (b_{10}) & 97\cdot 41\cdot 19\cdot 4 & (2,2,2,2,8)\\\hline
\rsp 64.150 &  (c_3)  & 3\cdot 8\cdot 11\cdot 23  &  (2,4,4) \\
\rsp       &  (d_1)  & 23\cdot 7\cdot 47\cdot 4  &  (2,8,8) \\
\end{array}$$
\end{table}\end{tiny}
\bigskip

\section{{\bf On the Length of the Narrow $2$-Class Field Tower of ${k}$}}

We are now interested in determining the length of the narrow $2$-class field tower of the real quadratic fields $k$ that we have been studying; and we have actually made some progress already. Namely, for all those $k$ for which $G^+/G^+_3$ are of order $32$, but not of type $32.033$, the length of the tower is $2$, by Corollary~\ref{C1} and Table 1 above, since $G^+/G^+_2$ is elementary and $G^+_2$ has rank $2$, which implies that $G^+$ is metabelian (refer also to the Main Theorem in \cite{BS2}).

On the other hand, when $G^+_2$ has rank at least $3$, we now wish to determine when  the length of the tower is $2$ and when it is at least $3$. This is the purpose of Part II of this project.

\newpage

{\bf Addresses}

\medskip

\author{Elliot Benjamin}

\address{School of Social and Behavioral Sciences,  Capella University,

  Minneapolis, MN 55402, USA; www.capella.edu}

\email{ben496@prexar.com}

\medskip

\author{C.\,Snyder}

\address{Department of Mathematics and Statistics, UMaine, Orono, ME 04469, USA}

\email{wsnyder@maine.edu}


\bigskip


\begin{thebibliography}{9999}
\bibitem{BLS0} E.\,Benjamin, F.\,Lemmermeyer, C.\,Snyder, {\em Imaginary Quadratic Fields $k$ with Cyclic
$\Cl_2(k^1)$}, J. Number Theory, {\bf 67} (1997), 229-245.
\bibitem{BLS} E.\,Benjamin, F.\,Lemmermeyer, C.\,Snyder, {\em Real Quadratic Fields with Abelian $2$-Class Field Tower},
J. Number Theory, {\bf 73} (1998), 182-194.
\bibitem{BLS2} E.\,Benjamin, F.\,Lemmermeyer, C.\,Snyder, {\em Imaginary Quadratic Fields with $\Cl_2(k)\simeq
(2,2,2)$}, J. Number Theory, {\bf 103} (2003), 38-70.
\bibitem{BS} E.\,Benjamin, C.\,Snyder, {\em Real Quadratic Fields with $2$-Class Group of Type $(2,2)$}, Math. Scan. {\bf 76} (1995), 161-178.
\bibitem{BS0} E.\,Benjamin, C.\,Snyder, {\em Some Real Quadratic Number Fields whose $2$-Class Fields have Class Number Congruent to $2$ Modulo $4$,} Acta Arith., {\bf 177} (2017), 375-392.
\bibitem{BS1} E.\,Benjamin, C.\,Snyder, {\em The Narrow $2$-Class Field Tower of some Real Quadratic Number Fields},  Acta Arith., {\bf 194.2} (2020), 187-218.
\bibitem{BS2} E.\,Benjamin, C.\,Snyder, {\em On some  Finite 2-Groups  in which the Derived Group has two Generators}, Czechoslovak Mathematical Journal, {\bf 73} (148) (2023), 71-100.
\bibitem{Bla2} N.\,Blackburn, {\em On Prime-Power Groups in which the Derived Group has Two Generators,} Proc.\,Camb.\,Phil.\,Soc., {\bf 53} (1957), 19-27.
\bibitem{Bla} N.\,Blackburn, {\em On Prime-Power Groups with two Generators}, Proc.\,Camb.\,Phil.\,Soc., {\bf 54}  (1958), 327-337.
\bibitem{CD} R.\,Couture, A.\,Derhem, {Un probl\`{e}m de capitulation}, C.\,R.\,Acad.\,Sci.\,Paris, {\bf 314}, S\'{e}rie I (1992), 785-788.
\bibitem{HS} M.\,Hall, J.\,K.\,Senior, {\em The Groups of Order $2^n$ $(n\leq 6)$}, Macmillan, New York, 1664.
\bibitem{Hup} B.\,Huppert, {\em Endliche Gruppen I}, Springer-Verlag, Heidelberg, 1967.
\bibitem{Kis} H.\,Kisilevsky, {\em Number fields with class number congruent to $4$ mod $8$ and Hilbert's Theorem $94$,} J. Number Theory, {\bf 8} (1976), 271-279.
\bibitem{Koc} H.\,Koch, {\em \"Uber den $2$-Klassenk\"orperturm eines quadratischen Zahlk\"orpers}, J.\,Reine Angew.\,Math., {\bf 214/215} (1963), 201-206.
\bibitem{Kub} T.\,Kubota, {\em \"Uber den Bizyklisch Biquadratischen Zahlk\"orper,} Nagoya Math.\ J.\ {\bf 10} (1956), 65-85.
%\bibitem{Lan} S.\,Lang, {\em Algebraic Number Theory}, Addison-Wesley, Reading, MA (1970).
\bibitem{Lem2} F.\,Lemmermeyer, {\em Kuroda's class number formula,} Acta Arith., {\bf 66.3} (1994), 245-260.
\bibitem{Lem} F.\ Lemmermeyer, {\em Die Konstruktion von Klassenk\"orpern,} Diss.\,Univ.\,Heidelberg (1995).
\bibitem{Lem3} F.\,Lemmermeyer, {\em Ideal Class Groups of Cyclotomic Number Fields}, Acta Arith., {\bf 72.4} (1995), 347-359.
\bibitem{Lem1} F.\ Lemmermeyer, {\em Reciprocity Laws}, Springer Verlag, 2000.
\bibitem{RR} L. R\'edei, H. Reichardt, {\em Die Anzahl der durch $4$ teilbaren Invarianten der Klassengruppe
    eines beliebigen quadratischen Zahlk\"orpers},  J. Reine Angew. Math. {bf 170} (1933), 69--74.
\bibitem{SW} T.\,W.\,Sag, J.\,W.\,Wamsley, {\em Minimal Presentations for Groups of Order $2^n$, $n\leq 6$}, J.\,Austral.\,Math.\,Soc., {\bf 15.4} (1973), 461-469.
\end{thebibliography}
\end{document}